\documentclass[3p, 11pt]{elsarticle}

\makeatletter
\def\ps@pprintTitle{
 \let\@oddhead\@empty
 \let\@evenhead\@empty
 \def\@oddfoot{}
 \let\@evenfoot\@oddfoot}
\makeatother

\usepackage{lipsum}
\usepackage{hyperref}
\usepackage{amssymb}
\usepackage{caption,subcaption}
\usepackage{amsmath}
\usepackage{pgfplots}
\usepackage{pgf}
\pgfplotsset{compat=1.18}
\usepackage{tikz} 
\usepackage{amsfonts}
\usepackage{graphicx}
\usepackage{epstopdf}
\usepackage{MnSymbol}
\usepackage{overpic}
\usepackage{float}
\usepackage{colortbl}
\usepackage[absolute,overlay]{textpos}
\usepackage{booktabs}
\usepackage{relsize}
\usepackage{mathrsfs}
\usepackage{algpseudocode}
\usepackage{mathtools}

\usepackage[linesnumbered, ruled, vlined]{algorithm2e}
\biboptions{sort&compress}
\usepackage{cleveref}
\usepackage{url}

\newcommand{\bm}[1]{\ensuremath{\mathbf{#1}}}

\newtheorem{definition}{Definition}

\journal{Computer Methods in Applied Mechanics and Engineering}

\begin{document}

\begin{frontmatter}

\title{A CUR Krylov Solver for Large-Scale Linear Matrix  Equations}

\author{Saeed Akbari$^{1}$, Damiano Lombardi$^{2}$, and Hessam Babaee$^{1*}$}

\address{$^{1}$ Department of Mechanical Engineering and Materials Science, University of
Pittsburgh, 3700 O’Hara Street, Pittsburgh, PA, 15213, USA \\ \vspace{2mm}
$^{2}$ COMMEDIA, Laboratoire Jacques–Louis Lions, Sorbonne Université et Inria Paris, 2 rue Simone Iff, 75012 Paris, France \\ \vspace{2mm}
* Corresponding Author, Email: h.babaee@pitt.edu \vspace{-8mm}}

\begin{abstract}
Developing efficient solvers for large-scale multi-term linear matrix equations remains a central challenge in numerical linear algebra and is still largely unresolved. This paper introduces a methodology leveraging CUR decomposition for solving large-scale generalized Sylvester as well as non-Sylvester multi-term equations on low-rank matrix manifolds. The approach decomposes the original equation into two smaller subproblems: one involving all columns with a small subset of rows, and the other involving all rows with a small subset of columns. The rows and columns are strategically selected using the discrete empirical interpolation method. We further utilize the CUR properties and propose a novel iterative scheme that removes the dependencies between selected and unselected rows (and likewise for columns), thereby enabling the subset problems to be solved independently.  We present a Krylov-based scheme for solving the resulting subproblems, which scales effectively to large problems and does not rely on a Sylvester structure.  The method incorporates rank adaptivity, dynamically adjusting computational rank to reach the desired accuracy.  The methodology is demonstrated in three representative settings: (i) implicit time integration of matrix differential equations on low-rank manifolds, leading to multi-term linear matrix equations; (ii) large-scale steady-state generalized Lyapunov equations including cases of size up to $10^{13}$ unknown entries; and (iii) non-Sylvester linear matrix equations with Hadamard product terms, such as those arising in nonlinear partial differential equations. 
\end{abstract}


\end{frontmatter}

\section{Introduction} \label{Introduction}
General multiterm linear matrix equations (LMEs) play a central role in a wide range of scientific and engineering applications, particularly in the discretization of deterministic and stochastic partial differential equations (PDEs) ~\cite{baumann2018msss,powell2017efficient}, PDE-constrained optimization~\cite{stoll2015low}, data assimilation~\cite{freitag2018low} and signal processing \cite{simoncini2016computational}.
Among these, Sylvester and Lyapunov equations have attracted particular attention because of their role in control theory and singnal processing ~\cite{benner2004solving,gajic2008lyapunov,laub1985numerical,van1991structured}, iterative linearization procedures for solving nonlinear matrix equations, such as algebraic Riccati equations~\cite{bini2011numerical}, and in the discretization of elliptic and parabolic PDEs~\cite{palitta2021matrix,palitta2016matrix,boulle2020computing,fortunato2020fast}. Sylvester equations also appear in applications such as eigenvalue placement for vibrating structures~\cite{brahma2009optimization}  and image denoising~\cite{calvetti1996application}.

The generalized Sylvester equation is given by 
\begin{equation}\label{eq:gse}
\bm A_1 \bm X \bm B_1 + \bm A_2 \bm X \bm B_2  = \bm C,   
\end{equation}
where \( \bm X \in \mathbb{R}^{n_1\times n_2} \) is the unknown matrix, and \( \bm A_i \in \mathbb{R}^{n_1 \times n_1} \), \( \bm B_i \in \mathbb{R}^{n_2 \times n_2} \) (for \( i=1,2 \)), and \( \bm C \in \mathbb{R}^{n_1 \times n_2} \) are known coefficient matrices. The generalized Sylvester equation encompasses several well-known matrix equations as special cases. Specifically, setting \( \mathbf{A}_2 = \mathbf{B}_1 = \mathbf{I} \) and \( \mathbf{A}_1 = \mathbf{A}, \; \mathbf{B}_2 = \mathbf{B} \) recovers the classical Sylvester equation~\cite{sylvester1884equation}
\begin{equation}\label{eq:stSylEq}
\bm{A} \bm{X} + \bm{X} \bm{B}  = \bm{C}.
\end{equation}
On the other hand, by choosing \( \mathbf{A}_1 = \mathbf{A} \), \( \mathbf{B}_1 = \mathbf{B} \), \( \mathbf{A}_2 = \mathbf{B}^{\mathrm T} \), and \( \mathbf{B}_2 = \mathbf{A}^{\mathrm T} \), Eq.~\eqref{eq:gse} reduces to the generalized Lyapunov equation
\begin{equation}\label{eq:gle}
\bm A \bm X \bm B + \bm B^{\mathrm T} \bm X \bm A^{\mathrm T}  = \bm C,
\end{equation}
where $\bm C$ is a symmetric negative definite matrix. This equation simplifies to the standard Lyapunov form when \( \mathbf{B} = \mathbf{I} \) as follows:
\begin{equation}\label{eq:stLyaEq}
\bm A \bm X + \bm X \bm A^{\mathrm T}  = \bm C.\end{equation}
Lyapunov equations have been studied extensively because of their central role in control theory, and as a result, numerical methods for their solution have progressed more rapidly than those for the general Sylvester equation.
Nonetheless, many of these techniques can be adapted to the Sylvester equation with some modifications. A detailed survey of Lyapunov equation techniques and their connections to control applications can be found in~\cite{gajic2008lyapunov}. The generalized Lyapunov Eq.~\eqref{eq:gle} can be transformed into the standard form
\begin{equation}\label{eq:ngle}
\tilde{\bm A} \bm X + \bm X \tilde{\bm A}^{\mathrm T}  = \tilde{\bm C},
\end{equation}
through left and right multiplication by \( \mathbf{B}^{-1} \) and \( \mathbf{B}^{\mathrm -T} \), respectively.
Regarding the existence and uniqueness of solutions, it is known from Sylvester equation theory that Eq.~\eqref{eq:ngle} admits a unique solution if and only if \( \lambda_i + \lambda_j \neq 0 \) for all eigenvalues \( \lambda_i, \lambda_j \) of \( \tilde{\bm A} \)~\cite{horn1991topics}.
This condition is naturally satisfied whenever 
\(\tilde{\bm A} \) is stable, that is, when all of its eigenvalues have negative real parts, which is a property commonly encountered in control applications.
 
For small-scale Sylvester and Lyapunov equations, typically involving matrices up to a few thousand in dimension, direct solvers are generally most efficient. A widely used method is the Bartels–Stewart algorithm~\cite{bartels1972solution}, which applies Schur decompositions to reduce the equations to triangular (or quasi-triangular) form, enabling efficient direct solutions. This algorithm is implemented in standard numerical libraries such as LAPACK and SLICOT~\cite{benner2009slicot,slowik2007evaluation,van1991slicot}.
For Lyapunov equations, additional speedups can be achieved by exploiting symmetry, which effectively reduces the computational cost by half~\cite{miller2000differential}. Hammarling’s method computes a Cholesky factor \(\mathbf{X} = \mathbf{L}\mathbf{L}^\top\) directly, where \(\mathbf{L}\) is a lower triangular matrix, avoiding explicit formation of the solution and improving numerical stability, especially for ill-conditioned cases~\cite{hammarling1982numerical}. A comparison of Hammarling’s and Bartels–Stewart’s methods is provided in~\cite{zhou2002numerical}, and further performance gains have been made via block versions for discrete-time Lyapunov equations~\cite{kressner2008block}.

In cases where the Sylvester equation involves matrices with one dimension relatively small (typically fewer than 1000 rows or columns) and the other very large, the solution strategy often combines direct methods along the smaller dimension with iterative methods along the larger dimension. Performing a Schur decomposition on the large-scale matrix is computationally prohibitive due to dense operations, requiring \(\mathcal{O}(n^3)\) time and \(\mathcal{O}(n^2)\) memory. Instead, one can apply a direct decomposition method, such as Schur decomposition to the smaller side, while using iterative, projection-based techniques for the larger side.

Solving large-scale Sylvester equations, where both the number of rows and columns are large, is considerably more challenging than cases that are large in only one dimension, primarily due to memory demands. While the coefficient matrices are often sparse, the solution matrix is typically dense. Direct solvers, such as the Bartels–Stewart algorithm, become impractical in this regime—not only because of their cubic floating-point operation (flop) cost, but also because storing the dense solution matrix $\bm{X}$ can exceed available memory. For example, solving a Lyapunov equation for a dynamical system with $n=10^6$ requires storing $10^{12}$ entries. Consequently, most successful large-scale approaches rely on low-rank approximations of $\bm{X}$ to alleviate memory requirements, coupled with iterative methods to reduce the flops cost.

For large-scale Sylvester equations, Krylov-based iterative solvers are typically used. 
Early approaches projected the residual onto polynomial Krylov spaces~\cite{saad1990numerical, jaimoukha1994krylov}, which often suffer from slow convergence. To address this limitation, extended Krylov methods were developed~\cite{simoncini2007new}, incorporating both matrix and inverse powers to accelerate convergence.
In~\cite{el2024krylov}, a high-order adaptive-rank implicit integrators with an extended Krylov solver was introduced for stiﬀ time-dependent PDEs with implicit schemes. Another direction of improvement is offered by rational Krylov subspaces, first proposed by Ruhe~\cite{ruhe1984rational} for eigenvalue approximation, where carefully chosen shifts were shown to accelerate convergence toward targeted spectral regions. In the setting of matrix functions and matrix equations~\cite{druskin1998extended,guttel2013black,knizhnerman2010new,guttel2013rational}, rational Krylov spaces were also reported to achieve faster convergence~\cite{beckermann2021low}.
When shift parameters are selected adaptively~\cite{druskin2011adaptive}, rational Krylov methods can outperform the extended Krylov approach with much smaller subspace size.
Restarting techniques~\cite{saad1986gmres,frommer1998restarted, simoncini2003restarted, frommer2017block} and problem-specific preconditioning~\cite{anzt2016updating, bellavia2011efficient} have also been introduced to manage memory and improve convergence.

Alternating Direction Implicit (ADI) iterations are widely used for solving large-scale Sylvester and Lyapunov equations~\cite{benner2008numerical, gugercin2003modified, li2002low, lu1991solution, penzl1999cyclic, wachspress1988iterative, benner2009adi}. Originating from classical operator-splitting schemes~\cite{peaceman1955numerical, ellner1986new}, ADI methods iteratively construct a low-rank factor that approximates the solution matrix. In their low-rank form, they can achieve performance comparable to rational Krylov subspace methods when good shift parameters are chosen—whether based on approximation theory or spectral heuristics. However, selecting robust shift parameters for general problems remains  challenging~\cite{benner2014computing}.

For a large-scale Sylvester equations, low-rank approximations are particularly attractive as they can reduce the storage requirements as well as flop costs.
For example, when the right-hand side \( \bm C \) in the Lyapunov equation is low rank, expressed as \( \bm C = -\bm B \bm B^T \) based on Eq.~\eqref{eq:lyaCont} with \( \bm B \in \mathbb{R}^{n \times n_u} \) and \( n_u \ll n \), the solution \( \bm X \) can often be approximated by a low-rank factorization \( \bm X \approx \bm Z \bm Z^T \), where \( \bm Z \in \mathbb{R}^{n \times r} \), \( r \ll n \), thereby requiring storage only for the smaller matrix \( \bm Z \)~\cite{antoulas2002decay, benner2013low, grasedyck2004existence}.
This idea has inspired the development of several iterative schemes~\cite{simoncini2007new} and low-rank ADI methods~\cite{benner2008numerical,li2002low,penzl2000eigenvalue}, which operate directly on the low-rank factors $\bm Z$ and make it possible to solve problems with dimensions as large as one million.

Beyond the classical Sylvester form, broader classes of LMEs arise in diverse applications—most notably in the implicit time integration of linear and nonlinear matrix differential equations (MDEs). Dynamical low-rank approximation (DLRA)~\cite{koch2007dynamical} provides a mathematical framework for solving such MDEs on a low-rank manifold. The method projects the right-hand side of the MDE onto the tangent space, thereby constraining the dynamics to remain on the low-rank manifold. However, it is well known that time integration of DLRA equations can become unstable in the presence of very small singular values. To address this issue, projector-splitting techniques~\cite{lubich2014projector} and, more recently, basis-update Galerkin (BUG) integrators~\cite{CL22} have been developed; both approaches remain robust when singular values are small or even zero. An alternative strategy is step truncation~\cite{RDV22}, which avoids tangent-space projections altogether by applying an singular value decomposition (SVD) to the full-order time-discretized solution, thereby preventing uncontrolled rank growth at each time step.

However, implicit time integration of nonlinear MDEs remains computationally expensive, since the nonlinear image of a low-rank matrix can be high-rank or even full rank. For instance, the element-wise exponential $\exp(\bm X)$ of a low-rank matrix $\bm X$ is generally full rank. To mitigate these costs, CUR decompositions have recently been demonstrated as an efficient and scalable alternative for explicit time integration of nonlinear MDEs~\cite{DPNFB23} and their extension to tensor differential equations~\cite{GBTT24}. CUR low-rank approximations, first introduced in~\cite{GTZ97} and also referred to as \emph{pseudoskeleton} or \emph{cross} approximations, differ from the SVD in their construction. Whereas the SVD expresses singular vectors as linear combinations of all rows and columns, CUR forms a rank-$r$ representation directly from only $\mathcal{O}(r)$ carefully selected rows and columns of the target matrix~\cite{MD09}. For a comprehensive survey of CUR methods, see~\cite{HH20}.

Recent work has introduced CUR-based time integration methods for implicit schemes applied to parametric PDEs on low-rank matrix manifolds~\cite{NAB25}. In this approach, linearized matrix equations are solved within Newton’s method. The LMEs obtained by matricizing discretized parametric PDEs are not of Sylvester type, as they involve only left-side matrix multiplication; hence, each column can be solved independently, a property exploited in~\cite{NAB25}. By contrast, non-parametric PDEs on low-rank manifolds give rise to equations with both left- and right-side multiplications, yielding multi-term LMEs. As we show in this work, such equations may also include Hadamard-product terms, distinguishing them from generalized Sylvester equations. Developing efficient solvers for these equations is the main objective of this work.

In this paper, we present a novel methodology for solving large-scale generalized Sylvester equations, multi-term Sylvester equations, and more general non-Sylvester LMEs. Our approach leverages CUR low-rank approximation, in which the LME is solved for a strategically selected set of subcolumns and subrows. This property, combined with a Krylov subspace solver, enables the efficient solution of large-scale LMEs. The main contributions of our work are:

\begin{itemize}
\item We present an iterative CUR low-rank approximation, which involves solving a time-dependent multi-term linear MDE for strategically selected rows and columns of the matrix. We presented a novel approach to eliminate the dependencies of the rows and columns on other rows and columns. 
\item We present a Krylov solver that exploits the Kronecker delta structure of the resulting row and column equations.
\item We extend the algorithm to address multi-term linear Sylvester equations as well as more general non-Sylvester forms. In particular, we demonstrate the solution of nonlinear PDEs using Newton’s method, where each iteration gives rise to a non-Sylvester multi-term linear equation.
\end{itemize}

The remainder of this paper is organized as follows. Section~\ref{sec:def} introduces notation and definitions. Section~\ref{sec:method} develops the methodology, covering the problem setup; optimal low-rank approximation via SVD; CUR and near-optimal CUR approximations; Krylov subspace solvers; and both Sylvester and non-Sylvester linear matrix equations. Section~\ref{sec:cases} presents demonstration cases, including linear and nonlinear transient heat equations and a steady-state Lyapunov equation. Upon discretization, these problems turn to generalized Sylvester, generalized Lyapunov, and non-Sylvester linear matrix equations, including a very large Lyapunov instance with approximately \(1.8\times10^{8}\) scalar algebraic equations. Section~\ref{sec:conc} concludes the paper.
\section{Notation and Definitions} \label{sec:def}
In this section, we introduce the notation used throughout the paper. Vectors are denoted by boldface lowercase letters, e.g., \( \bm{x} \in \mathbb{R}^n \), and matrices by boldface uppercase letters, e.g., \( \bm{X} \in \mathbb{R}^{n_1 \times n_2} \). We denote the identity matrix of size $n \times n$ by $\bm I_n$.

We use MATLAB indexing notation to extract a subset of rows and columns of a matrix. To explain this notation, let $\bm p =[p_1,p_2, \dots, p_r ]$ and $\bm q =[q_1,q_2, \dots, q_r ]$ be  vectors containing row and column indices, respectively. Then, for matrix $\bm X \in \mathbb{R}^{n_1 \times n_2}$, 
$\bm X(\bm p,:) \in \mathbb{R}^{r \times n_2}$ selects all columns at the $\bm p$ rows, and $\bm X(:,\bm q)\in \mathbb{R}^{n_1 \times r}$ selects all rows at the $\bm q$ columns of the matrix $\bm X$. We also use indexing matrices $\bm P = \bm I_{n_1}(\bm p,:) \in \mathbb{R}^{r \times n_1}$ and  $\bm Q = \bm I_{n_2}(\bm q,:) \in \mathbb{R}^{r \times n_2}$. It is easy to verify that: $\bm X(\bm p,:) = \bm P \bm X$ and $\bm X(:,\bm q) =  \bm X \bm Q$.  

In the following, we define the low-rank matrix manifolds. 
 \begin{definition}Low-rank Matrix Manifolds : \label{def:Mr}
Let $\mathcal{M}_r$ be the set of all $n_1 \times n_2$ real-valued matrices with rank equal to $r$, for some integer $r$. Then $\mathcal{M}_r$ is called a low-rank matrix manifold of rank $r$.
\begin{equation*}
\mathcal{M}_r = \{\hat{\bm X} \in \mathbb{R}^{n_1 \times n_2}: \ \mbox{rank} (\hat{\bm X}) = r \}.
\end{equation*}
The hat symbol $( \hat{ \ \ } )$ is used to represent rank-$r$ matrices belonging to the set $\mathcal{M}_r$.
\end{definition}

We define a vector \( \bm{x} \in \mathbb{R}^{n} \) as a column of \( n \) real-valued entries, i.e., an element of the Euclidean space \( \mathbb{R}^{n} \). A matrix \( \bm{X} \in \mathbb{R}^{n_1 \times n_2} \) is a two-dimensional array of real numbers with \( n_1 \) rows and \( n_2 \) columns.
We denote the Frobenius norm of a matrix with \( \| \cdot \|_F \).

\section{Methodology} \label{sec:method}

\subsection{Problem setup}
We consider the linear matrix differential equations as follows:
\begin{equation}\label{eq:gse_td}
\bm A_0(t) \dot{\bm  X}(t) \bm B_0(t)= \sum_{i=1}^d\bm A_i(t) \bm X (t) \bm B_i(t)  + \mathbf{E}(t) \circ \bm X(t) - \bm C(t),   
\end{equation}
given an initial condition given by $\bm X (0) = \bm X_0$, where $\bm X \in \mathbb{R}^{n_1\times n_2}$, $\bm A_i \in \mathbb{R}^{n_1 \times n_1}$ and $\bm B_i \in \mathbb{R}^{n_2 \times n_2}$ for $i=0, \dots, d$,  $\bm E \in \mathbb{R}^{n_1\times n_2}$, $\circ$ denotes the Hadamord (elementwise) product and $\dot{\bm X} = d\bm X/dt$. We consider the general case where the matrices $\bm A_i$, $\bm B_i$, $\bm E$, and $\bm C$ are time-dependent. Typically,  $\bm{A}_i$  and  $\bm{B}_i$  are sparse, while  $\bm{C}$  may be a non-sparse, full-rank matrix. 

The above equation arises from the discretization of time-dependent PDEs. Differential generalized Sylvester and Lyapunov equations can be seen as a special case of the above equation \cite{BBH19}. For example,  when $d=2$, and $\bm E \equiv \bm 0$, Eq.~\eqref{eq:gse_td} represents generalized differential Sylvester and forward Lyapunov equations when $\bm B_i = \bm A_i ^{\mathrm T}$.  The steady-state generalized Sylvester given by Eq.~\eqref{eq:gse} is a special case where $\bm A_i$, $\bm B_i$,  and $\bm C$ are time-invariant and as $t\rightarrow \infty$, the solution of Eq.~\eqref{eq:gse_td} approaches to the solution of Eq.~\eqref{eq:gse}.

Our approach approximates $\mathbf{X}(t)$ in a low-rank form, which amounts to integrating Eq.~\eqref{eq:gse_td} on the low-rank matrix manifold $\mathcal{M}_r$. We focus on implicit time integration of Eq.~\eqref{eq:gse_td} on  $\mathcal{M}_r$. This choice has two key advantages: first, implicit integration of MDEs on low-rank manifolds is an important problem in its own right as it appears in the time-integration of stiff PDEs or stiff differential Lyapunov equations; second, it enables stable time stepping with relatively large $\Delta t$, which is especially useful when seeking steady-state solutions.

 For the sake of simplicity in the exposition, we present our methodology for the case of $d=2$, and $\bm E  \equiv \bm 0$. This special case covers some of the most common LMEs, namely, the differential generalized Sylvester and Lyapunov equations as well as time MDEs arising from the spatial discretization of PDEs. The development for this special case involves all the elements of the methodology. Extension of the algorithm to $d>2$ and $\bm E \neq \bm 0$ is straightforward and is discussed in \S \ref{sec:non_gen_syl}. 
 
 For the sake of simplicity, we consider implicit Euler time advancement as follows:
\begin{equation}
\bm A_0^{k+1}\frac{(\bm X^{k+1} - \bm X^{k})}{\Delta t} \bm B_0^{k+1} = \bm A_1^{k+1}\bm X^{k+1}\bm B_1^{k+1} +  \bm A_2^{k+1}\bm X^{k+1}\bm B_2^{k+1} - \bm C^{k+1}, 
\end{equation}
where the superscript denotes the time step, e.g., $ \bm A_1^{k} = \bm A_1(t_k)$. Extension to higher-order implicit time integration is discussed later in the paper. For brevity in the notation, we omit the explicit dependence of the coefficient matrices on the time step, with the understanding that all are evaluated at $t_{k+1}$.  Rearranging the above equation results in:
\begin{equation}
\bm A_0 \bm X^{k+1} \bm B_0- \Delta t\big [\bm A_1\bm X^{k+1}\bm B_1 +  \bm A_2\bm X^{k+1}\bm B_2 \big ]=  \bm A_0 \bm X^{k}  \bm B_0- \Delta t  \bm C.
\label{eq:trgse}
\end{equation}
The above equation can be recast in the form of:
\begin{equation}
\mathcal A (\bm X^{k+1}) = \mathcal B (\bm X^{k}),
\end{equation}
where $\mathcal A (\bm X): \mathbb{R}^{n_1 \times n_2} \rightarrow \mathbb{R}^{n_1 \times n_2}$ and    $\mathcal B (\bm X): \mathbb{R}^{n_1 \times n_2} \rightarrow \mathbb{R}^{n_1 \times n_2}$  defined as:
\begin{subequations}
\begin{align}
    \mathcal A (\bm X) &=  \bm A_0\bm X \bm B_0 - \Delta t\big [\bm A_1\bm X\bm B_1 +  \bm A_2\bm X\bm B_2 \big ] \label{eq:LME_A}\\
    \mathcal B (\bm X)  &= \bm A_0 \bm X \bm B_0  - \Delta t  \bm C\label{eq:LME_B}
    \end{align}
\end{subequations}
The above equation is a multi-term LME. In many applications, such as finite-element discretizations of PDEs, one often has $\mathbf{A}_0 = \mathbf{A}_i$ or $\mathbf{B}_0 = \mathbf{B}_i$ for $i=1,2$, by absorbing the term $\mathbf{A}_0 \mathbf{X} \mathbf{B}_0$ into one of the $\mathbf{A}_i \mathbf{X} \mathbf{B}_i$ terms. In such cases, the LME can be written as a generalized Sylvester (two-term) equation. Here, however, we focus on the more general three-term case, since the proposed approach is agnostic to the number of terms.

For very large $n_1$ and $n_2$, i.e., $n_1,  n_2 \sim \mathcal{O}(10^4)$, or larger matrices, solving the above LME becomes cost-prohibitive. In such cases, even storing  $\bm{X}^{k}$  may be unfeasible since $\bm{X}^{k}$ is not a sparse matrix even when the coefficient matrices are sparse.

The high computational cost of solving the LME motivates the use of low-rank approximation, in which $\bm{X}^{k}$ is approximated by a rank-$r$ matrix, denoted by $\hat{\bm{X}}^k \in \mathbb{R}^{n_1 \times n_2}$, with $r \ll n_1,n_2$. A rank-$r$ matrix can be parameterized as the product of two tall matrices, requiring only $r(n_1 + n_2)-r^2$ parameters—significantly fewer than the original $n_1 n_2$ parameters. This effectively amounts to solving Eq.~\eqref{eq:gse_td} on a low-rank matrix manifold, $\mathcal{M}_r$.

Several methods exist for integrating such equations on low-rank manifolds, including the dynamical low-rank approximation framework via orthogonal projection onto the tangent space \cite{KL07}, the step truncation method using orthogonal projection \cite{RDV22}, and the oblique projection approach based on the CUR algorithm \cite{DPNFB23,D25}.

In this work, we adopt the step truncation approach with oblique projection via CUR, primarily for its computational efficiency, which will be discussed in detail in later sections, as well as the simplicity of its implementation. Our algorithm consists of two main components: CUR-based low-rank approximation and a Krylov subspace solver. To set the stage, we begin by presenting the optimal low-rank approximation via SVD. While SVD is not computationally efficient, it provides a baseline for optimality that informs our subsequent developments.

\subsection{CUR low-rank approximation}
We begin by presenting an efficient strategy for matrix approximation that builds upon a certain kind of CUR low-rank approximation. A CUR algorithm was first introduced in \cite{GTZ97} -also referred to as the pseudoskeleton or cross approximation.    This approach addresses the computational challenges associated with optimal low-rank approximation via SVD by selecting subsets of rows and columns to construct interpolatory low-rank approximations. The CUR (cross) algorithm is briefly explained below.

Let $\mathbf{X} \in \mathbb{R}^{m \times n}$ be a matrix, and let $\bm p = \{ i_\alpha^1 \}$ and $\bm q = \{ i_\alpha^2 \}$ (for $\alpha = 1, \ldots, r$) be index sets corresponding to selected rows and columns.
Let the column and row submatrices be denoted as  $\mathbf{\tilde C} = \mathbf{X}(:,\bm q)$ and $\mathbf{\tilde R} = \mathbf{X}(\bm p,:)$, respectively.
A CUR decomposition constructs a rank-$r$ approximation of a matrix $\mathbf{X} \approx \mathbf{\tilde C} \tilde{\mathbf{U}} \mathbf{\tilde R}$, where $\mathbf{\tilde C} \in \mathbb{R}^{n_1 \times r}$, $\mathbf{\tilde R} \in \mathbb{R}^{r \times n_2}$, and $\tilde{\mathbf{U}} \in \mathbb{R}^{r \times r}$. If the columns of \( \bm{\tilde  C} \) and the rows of \( \bm{\tilde R} \) are each linearly independent, then the CUR approximation yields a rank-\( r \) approximation. We denote the CUR algorithm as: 
\[
\hat{\bm X} = \texttt{CUR}(\bm X).
\]
The accuracy of the CUR low-rank approximation depends on the choice of the indices and how $\tilde{\mathbf{U}}$ is obtained. Given  the row and column selection, the best $\tilde{\mathbf{U}}$ (in the Frobenius norm) is obtained via the \emph{orthogonal projection} of the target matrix onto the space spanned by the selected columns and rows \cite{SE16}:
\begin{equation}
\tilde{\mathbf{U}}_{best} = (\mathbf{\tilde C}^\top \mathbf{\tilde C})^{-1} \mathbf{\tilde C}^\top \mathbf{X} \mathbf{\tilde R}^\top (\mathbf{\tilde R} \mathbf{\tilde R}^\top)^{-1}.
\end{equation}
However, this method requires full access to all entries of $\mathbf{X}$, which is undesirable when access to data is costly, for example, if we have to perform extra computation in obtaining the entries of $\bm X$. An alternative is to interpolate any column and row onto the row and column subspaces. This amounts to an \emph{oblique projection}, which requires access only to the intersection of the selected rows and columns of the original matrix $\mathbf{X}$, yielding:
\begin{equation}
\tilde{\mathbf{U}} = \mathbf{X}^{-1}(\bm p,\bm q).
\end{equation}
The resulting CUR approximation then becomes:
\begin{equation}
\hat{\mathbf{X}} = \mathbf{X}(:,\bm q) \mathbf{X}^{-1}(\bm p,\bm q) \mathbf{X}(\bm p,:). 
\label{eq:cur_interp}
\end{equation}
More details related to orthogonal versus oblique projection can be found in \cite{DPNFB23}.  This interpolatory CUR algorithm requires only $s = r(n_1 + n_2) - r^2$ matrix entries, which is the minimal number needed to construct a rank-$r$ approximation. For $n_1 = n_2=n$ and $r \ll n$, this count is roughly a factor of $n^2 / 2rn = n / 2r$ smaller than the full matrix size. It is easy to verify that the above CUR low-rank approximation satisfies the interpolation properties, i.e.,  $\mathbf{X}(\bm p,:) = \hat{\mathbf{X}}(\bm p,:)$ and $\mathbf{X}(:,\bm q) = \hat{\mathbf{X}}(:,\bm q)$. In other words, the CUR low-rank approximation is equal to the target matrix at the selected columns and rows. 

The selection of $\bm p$ and $\bm q$ could significantly affect the CUR approximation error, i.e., $\mathbf{X} - \hat{\mathbf{X}}$. We use the discrete empirical interpolation method (DEIM) \cite{CS10},    which has been successfully applied to solving nonlinear parametric PDEs on low-rank manifolds ~\cite{MNAdaptive, DPNFB23}.   There are also other excellent algorithms for column (and row) selections, including QDEIM \cite{DG16} and the max-vol algorithm \cite{GT01} that seeks to maximize the row-column intersection volume.  The DEIM algorithm chooses these indices by utilizing the exact or approximate dominant left and right singular vectors $\mathbf{U} \in \mathbb{R}^{n_1 \times r}$ and $\mathbf{Y} \in \mathbb{R}^{n_2 \times r}$ of $\mathbf{X}$, respectively. The DEIM operator is applied as $\bm p = \texttt{DEIM}(\mathbf{U})$ and $\bm q = \texttt{DEIM}(\mathbf{Y})$, and the overall method is referred to as CUR-DEIM. In iterative or time-dependent problems, the singular vectors are not available at the current iteration or time step. To address this, we use the singular vectors from the previous iteration to determine which rows and columns to sample. This approach naturally leads to a \emph{targeted adaptive sampling}, as the selected rows and columns evolve with each iteration or time step.

According to~\cite{DPNFB23}, the approximation error of CUR-DEIM is bounded by:
\begin{equation}
\| \mathbf{X} - \hat{\mathbf{X}} \|_2 \leq c \hat{\sigma}_{r+1},
\label{eq:curdeim_error}
\end{equation}
where $\hat{\sigma}_{r+1} = \max\{ \|(\mathbf{I} - \mathbf{U}\mathbf{U}^\top)\mathbf{X} \|_2, \| \mathbf{X} (\mathbf{I} - \mathbf{Y}\mathbf{Y}^\top) \|_2 \}$ represents the orthogonal projection error and $\| \cdot \|_2$ is the second matrix norm. The term $\hat{\sigma}_{r+1}$ approximates the $(r+1)$-th singular value of $\mathbf{X}$, providing a lower bound on the approximation error.
The constant $c \geq 1$ accounts for the condition number of the inverse of the interpolation matrices, given by:
\begin{equation}
c = \min\{ \eta_r (1 + \eta_c), \eta_c (1 + \eta_r) \},
\end{equation}
where $\eta_r = \| \mathbf{U}^{-1}(\bm p, :) \|_2$ and $\eta_c = \| \mathbf{Y}^{-1}(\bm q, :) \|_2$.  This  DEIM algorithm employs a greedy selection strategy for interpolation points that seeks to minimize $\eta_c$ or $\eta_r$. It is possible to perform oversampling to further reduce the norm of the inverse of the submatrices \cite{PB25}; however, no such oversampling strategies are employed in this work. 
 
\subsection{Optimal low-rank approximation based on SVD}
We consider a rank-$r$ approximation of Eq.~\eqref{eq:gse}. We denote the low-rank approximation of $\bm X$ with $\hat{\bm X}$.  Replacing $\hat{\bm X}$ in Eq.~\eqref{eq:gse} generates the residual as follows:
\begin{equation}
 R(\hat{\bm X}^{k+1}) =  \big \| \mathcal A (\hat{\bm X}^{k+1}) - \mathcal B (\hat{\bm X}^{k})  \big \|_F.
\end{equation}
The optimal solution to the above problem can be formulated as in the following: 
\[
\hat{\bm X}_{best}^{k+1}= \underset{\hat{\bm Y} \in \mathcal{M}_r}{\arg\min} \ R(\hat{\bm Y}).
\] 
The minimization problem outlined above is a constrained, nonconvex optimization problem. A brute-force approach to solving it involves first converting the matrix equation into vector form using the Kronecker product, as shown below:
\begin{align*}
\bm A &= \bm B_0^{\mathrm T} \otimes \bm A_0 - \Delta t\big [\bm B_1^{\mathrm T} \otimes \bm A_1 + \bm B_2^{\mathrm T} \otimes \bm A_2 \big ],\\
\bm b &= \mbox{vec}(\mathcal B (\hat{\bm X}^{k})),
\end{align*}
where $\bm A \in \mathbb{R}^{n_1 n_2 \times n_1 n_2}$ and $\bm b \in \mathbb{R}^{n_1 n_2 \times 1}$. Then, \( \bm{X}^{k+1} = \text{mat}(\bm{x}^{k+1}) \), where \( \bm{x}^{k+1} = \bm{A}^{-1} \bm{b} \). The best rank-\( r \) approximation, \( \hat{\bm{X}}_{\text{best}}^{k+1} \), can be obtained by computing the SVD of \( \bm{X}^{k+1} \) and truncating at rank \( r \).

Using the approach above is impractical, as it requires solving the full-order model (FOM) -- precisely what we aim to avoid. In the following, we present a near-optimal low-rank approximation based on the CUR-DEIM decomposition.

\subsection{Near-optimal low-rank approximation based on CUR}\label{subsec:nearOpt}
In the following, we present a CUR methodology to construct a near-optimal solution to the above residual minimization problem. Let us denote the residual matrix with:
\[
\bm R = \mathcal A (\hat{\bm X}^{k+1}) - \mathcal B (\hat{\bm X}^{k}).
\]
Our approach is a residual collocation method where the residual is set to zero at strategically selected columns and rows. Let $\bm p = [p_1, p_2, \dots, p_r]$ and $\bm q = [q_1, q_2, \dots, q_r]$ be the indices of $r$ rows and columns of matrix $\bm R$, respectively, and $r$ is the rank of $\hat{\bm X}$. Our approach constructs a CUR low-rank approximation for $\hat{\bm X}^{k+1}$. This requires solving for $\bm X^{k+1}(:,\bm q)$ and $\bm X^{k+1}(\bm p,:)$. We first show how $\bm X^{k+1}(:,\bm q)$ can be computed. We set the residual to zero at the DEIM-selected indices as follows:
\begin{equation}\label{eq:col}
\bm A_0 \bm X(:,\bm q) \bm B_0- \Delta t\big [\bm A_1\bm X\bm B_1(:,\bm q) +  \bm A_2\bm X\bm B_2(:,\bm q) \big ]=  \bm A_0 \bm X^{k}(:,\bm q) \bm B_0  - \Delta t  \bm C(:,\bm q).
\end{equation}
where for the sake of simplicity in notation we use $\bm X \equiv \bm X^{k+1}$.  The main challenge in solving the equation above lies in the fact that, as long as matrices \( \bm{B}_i \), $i=0,1,2$  are not diagonal, solving for \( \bm{X}(:, \bm{q}) \) depends on other columns of \( \bm{X} \). Ideally, we seek to derive a system of equations for \( \bm{X}(:, \bm{q}) \) alone. In the following, we present an iterative algorithm that resolves this dependency and enables solving \( \bm{X}(:, \bm{q}) \). When \( \bm{B}_i \), $i=0,1,2$ or alternatively, \( \bm{A}_i \), $i=0,1,2$ are diagonal, then column (or rows) can be solved independently, resulting in a simpler algorithm that was recently introduced \cite{NAB25}. 

Let us denote the nonzero entries of submatrix $\bm B_i(:,\bm q)$ with $\bm B_i(\bm q_{\bm B_i},\bm q)$ where $\bm q_{\bm B_i}$ is a vector of size $n_{\bm B_i}$ containing the union of the  indices at which $\bm B_i(:,\bm q)$ is nonzero. Therefore, 
\[
\bm X \bm B_i(:,\bm q) = \bm X(:,\bm q_{\bm B_i}) \bm B_i(\bm q_{\bm B_i},\bm q),
\]
We use the CUR low-rank approximation to express $\bm X(:,\bm q_{\bm B_i})$ versus $\bm X(:,\bm q) $. We use the fact that in CUR low-rank approximation, any column of the matrix is approximated in the span of the selected columns, i.e.,  $\bm X(:,\bm q) $. Therefore, the following relationship holds:
\begin{equation}\label{eq:col_aux1}
\bm X(:,\bm q_{\bm B_i}) \approx  \hat{\bm X}(:,\bm q_{\bm B_i}) = \bm X(:,\bm q) \bm Z_{\bm B_i},
\end{equation}
where $\bm Z_{\bm B_i} \in \mathbb{R}^{r\times  n_{\bm B_i}}$. However,  $\bm Z_{\bm B_i}$ is unknown and must be computed. It is possible to obtain $\bm Z_{\bm B_i}$ via 
\[
\bm Z_{\bm B_i} = \bm X(:,\bm q)^{\dagger} \bm X(:,\bm q_{\bm B_i}).
\]
However, the  $\bm X(:,\bm q)$ can become ill-conditioned when $\sigma_r$, i.e., the $r$th  singular value of $\bm X$ is very small. It is possible, however, to obtain a stable approximation of $\bm Z_{\bm B_i}$ in the presence of small or even zero singular values of $\bm X$. 
Replacing  $\bm X$ with its low-rank approximation in the above equation results in:
\[
\bm U \bm \Sigma \bm Y(\bm q_{\bm B_i},:)^{\mathrm T} = \bm U \bm \Sigma \bm Y(\bm q,:)^{\mathrm T}\bm Z_{\bm B_i}.
\]
Therefore, 
\[
\bm Z_{\bm B_i} =  [\bm Y^{\mathrm T}(\bm q,:)]^{-1} \bm Y^{\mathrm T}(\bm q_{\bm B_i},:).
\]
DEIM ensures that the $\bm Y(\bm q,:)]$ is invertible. In fact, DEIM is a greedy algorithm that seeks to minimize $\eta_c = \| \bm Y^{-1}(\bm q,:)\|$. 
Replacing Eq.~\eqref{eq:col_aux1} into Eq.~\eqref{eq:col} results in:
\begin{equation}\label{eq:col_aux3}
\bm A_0 \bm X(:,\bm q) \bm B_{0_r} - \Delta t\big [\bm A_1\bm X(:,\bm q)\bm B_{1_r} +  \bm A_2\bm X(:,\bm q) \bm B_{2_r} \big ]=  \bm A_0 \bm X^{k}(:,\bm q_{\bm B_0}) \bm B_0(\bm q_{\bm B_0}, \bm q)  - \Delta t  \bm C(:,\bm q),
\end{equation}
where  $\bm B_{i_r} \in \mathbb{R}^{r \times r}$ as follows: 
\[
\bm B_{i_r} = \bm Z_{\bm B_i} \bm B_i(\bm q_{\bm B_i},\bm q) \quad i=0,1,2.
 \]
 Computing the rows follows steps analogous to those for the columns, as shown below:
 \begin{equation}\label{eq:row}
\bm A_{0_r} \bm X(\bm p,:) \bm B_0 - \Delta t\big [\bm A_{1_r}\bm X(\bm p,:)\bm B_1 +  \bm A_{2_r}\bm X(\bm p,:) \bm B_2 \big ]=  \bm A_0 (\bm p,\bm p_{\bm A_0}) \bm X^{k}(\bm p_{\bm A_0},:)\bm B_0  - \Delta t  \bm C(\bm p,:),
\end{equation}
where $\bm A_{i_r} \in \mathbb{R}^{r \times r}$ are computed as follows: 
\[
\bm A_{i_r} =  \bm A_i(\bm p,\bm p_{\bm A_i}) \bm Z_{\bm A_i} \quad i=0,1,2.
 \]
 The matrix  $\bm Z_{\bm A_i}$ is determined such that:
 \[
 \bm X(\bm p_{\bm A_i},:) = \bm Z_{\bm A_i} \bm X(\bm p,:).
 \]
 Using the low-rank approximation results in:
 \[
\bm U(\bm p_{\bm A_i},:) \bm \Sigma \bm Y^{\mathrm T} = \bm Z_{\bm A_i} \bm U(\bm p,:) \bm \Sigma \bm Y^{\mathrm T}.
\]
From the above relationship, 
\[
\bm Z_{\bm A_i} = \bm U(\bm p_{\bm A_i},:) [\bm U(\bm p,:)]^{-1}.
\]
 To summarize, $\bm Z_{\bm A_i}$ and $\bm Z_{\bm B_i}$ are matrices are obtained from:
\begin{subequations}
\begin{align}
\bm Z_{\bm A_i} &= \bm U(\bm p_{\bm A_i},:) [\bm U(\bm p,:)]^{-1}, \label{eqn:ZA1}\\  
\bm Z_{\bm B_i} &=  [\bm Y^{\mathrm T}(\bm q,:)]^{-1} \bm Y^{\mathrm T}(\bm q_{\bm B_i},:) \label{eqn:ZB1}
.\end{align}
\end{subequations}
The matrices \( \bm{Z}_{\bm{A}_i} \)  and \( \bm{Z}_{\bm{B}_i} \) depend on the unknown matrices \( \bm{U} \) and \( \bm{Y} \). To address this, we propose a fixed-point iterative algorithm to determine these matrices, along with \( \bm{X}(:, \bm{q}) \) and \( \bm{X}(\bm{p}, :) \). We assume the previous time step is already stored in the low-rank form, i.e., $\hat{\bm X}^{k} = \bm U^k \bm \Sigma^k \bm Y^{k^\top}$. We denote the iteration count by the subscript $m$, i.e., $\hat{\bm X}^{k+1}_m$, and continue the iterations until convergence. The convergence criterion is described later in this section. Once the solution at time step $k+1$ satisfies this criterion, the final iterate is accepted as $\hat{\bm X}^{k+1}$, with the iteration subscript dropped.  

In summary, the coefficient matrices are updated according to the low-rank solution obtained from the previous iteration, as outlined below:
\begin{enumerate}
\item Compute $\bm Z_{\bm A_i}$  and $\bm Z_{\bm B_i}$ according to Eqs.~\eqref{eqn:ZA1}-\eqref{eqn:ZB1}. 
\item Solve Eqs.~\eqref{eq:col_aux3} and \eqref{eq:row} to obtain $\bm X(:,\bm q)$ and  $\bm X(\bm p,:)$, respectively.  
\item Use the sampled $\bm X(:,\bm q)$ and  $\bm X(\bm p,:)$ and update $\hat{\bm X}^{k+1}_m = \bm U \bm \Sigma \bm Y^{\mathrm T}$ using the CUR algorithm.
\item Update the row \( \bm p\) and column \( \bm q\) indices using the DEIM algorithm using the updated $\bm U$ and $\bm Y$.
\item $m \leftarrow m + 1$.
\item Repeat Steps 2-6 until convergence is achieved.
\end{enumerate}

In the above steps, the matrices $\bm Z_{\bm A_i}$, $\bm Z_{\bm B_i}$, $\bm U$, $\bm \Sigma$, and $\bm Y$ (representing $\hat{\bm X}^{k+1}_m$ in the compressed form), together with the DEIM indices $\bm p$ and $\bm q$, all vary with each iteration; for brevity, the iteration index $m$ and time-step index $k+1$ are omitted.

For the convergence criterion of $\hat{\bm X}^{k+1}$, we track the difference between successive iterates,  
\[
\Delta \bm X = \hat{\bm X}_{m+1}^{k+1} - \hat{\bm X}_m^{k+1}
\]
to assess the convergence of the inner CUR iterations. This ensures the convergence of $\bm Z_{\bm A_i}$ and $\bm Z_{\bm B_i}$ while simultaneously identifying suitable indices for $\bm p$ and $\bm q$. Numerical experiments reveal that convergence is reached quickly, typically within 5–15 iterations across a wide spectrum of problems.

To ensure that the LME is solved up to the desired accuracy, we monitor the residual due to low-rank approximation, 
\[
\bm R = \bm A_0 \hat{\bm X}^{k+1} \bm B_0- \Delta t\big [\bm A_1\hat{\bm X}^{k+1}\bm B_1 +  \bm A_2\hat{\bm X}^{k+1}\bm B_2 \big ]-  \bm A_0 \hat{\bm X}^{k} \bm B_0  + \Delta t  \bm C.
\]
If $\|\bm R \|_F$ exceeds the desired threshold, the rank is increased. 

In the following we propose an efficient algorithm to compute $\| \bm R\|_F$ and $\| \Delta \hat{\bm X}\|_F$. 
First, we note that explicitly forming the full residual matrix $\bm R \in \mathbb{R}^{n_1 \times n_2}$ and the update $\Delta \bm X \in \mathbb{R}^{n_1 \times n_2}$ is both unnecessary and memory-prohibitive. However, since $\hat{\bm X}^{k+1}$ is a rank-$r$ matrix, $\bm A_i \hat{\bm X}^{k+1} \bm B_i$ matrices are also  rank-$r$. Therefore, it is possible to add these rank-$r$ matrices, as well as a rank-$r$ approximation of $\mathbf C$, in the compressed form and obtain a rank-$5r$ matrix. However, there are cases that computing the residual in the low-rank form may become very costly. We consider two such scenarios: 
\begin{enumerate}
\item If the LME has a large number of terms, i.e., large $d$, the matrix of $\sum_{i=1}^d \bm A_i \hat{\bm X}^{k+1} \bm B_i$ is of rank $dr$. 
\item The matrix $\bm C$ may not be an exactly low-rank matrix, for example, it may be a time-dependent full-rank matrix, and the computation of its low-rank approximation may be costly.  
\end{enumerate}

To construct the low-rank approximation of \( \bm{R} \) and $\Delta\bm X$, we use CUR to evaluate these matrices at a small set of rows and columns. The number of sampled indices is equal to the estimated rank of the $\mathbf R$ and $\Delta\bm X$, denoted by \( r_{\text{res}} \) and \( r_{\Delta} \). The row and column indices are selected using  DEIM, applied to the singular vectors of the residual and $\Delta\bm X$ obtained in the previous iteration.

Once the required entries are computed, we compute  the CUR low-rank approximation of these matrices using \Cref{alg:SCUR} to construct a low-rank approximation in the SVD form:
\begin{equation}
\bm{R} \approx \hat{\bm R} = \bm{U}_R \bm{\Sigma}_R \bm{Y}_R^{\top},
\qquad
\Delta\bm X \approx \widehat{\Delta\bm X} = \bm{U}_\Delta \bm{\Sigma}_\Delta \bm{Y}_\Delta^{\top}. 
\end{equation}

Finally, we compute the Frobenius norm of them using the singular values \( \bm{\Sigma}_R \) and \( \bm{\Sigma}_\Delta \), avoiding the need to evaluate the full residual matrix:
\begin{equation}
\|\hat{\bm{R}}\|_F = \|\bm{\Sigma}_{R}\|_F = \sqrt{\sum_{i=1}^{r_{\text{res}}} \sigma_{R_i}^2},
\qquad
\|\widehat{\Delta\bm X}\|_F = \|\bm{\Sigma}_\Delta\|_F = \sqrt{\sum_{i=1}^{r_{\Delta}} \sigma_{\Delta_i}^2}.
\end{equation}

This procedure allows for efficient and scalable evaluation of residual norms in large-scale problems.

\Cref{alg:ITDB} summarizes the methodology presented in this work, which we refer to as the \texttt{TDB-CUR} algorithm, where TDB refers to the time-dependent (or adaptive) bases for the column and row spaces of $\bm X$.  
\begin{algorithm}[!h]
\begingroup
\fontsize{9pt}{9pt}\selectfont
\SetAlgoLined
\KwIn{$r$, $\epsilon_\Delta$, $\epsilon_R$, $\Delta r$, $\Delta r_R$, $\Delta r_{\Delta}$}
\KwOut{$\mathbf{U}^{k}$, $\boldsymbol \Sigma^{k}$, $\mathbf{Y}^{k}$}
Initialize $\mathbf{X}_{0} \in \mathbb{R}^{n_1 \times n_2}$\;
$\bm{U}^{k-1}$, $\bm \Sigma^{k-1}$, $\bm{Y}^{k-1}$  \hspace{7mm} $\rhd$ Taking SVD of the initial guess\;
\While{$\|\mathbf{R}\|_F > \epsilon_{R}$}{
    $\|\Delta\bm X\|_F \gets 1$\;
    \While{$\|\Delta\bm X\|_F > \epsilon_\Delta$}{
        $\bm{q} \gets \texttt{DEIM}(\mathbf{Y}^{k-1})$ \hspace{7mm} $\rhd$ column indices selection via DEIM\;
        $\bm{p} \gets \texttt{DEIM}(\mathbf{U}^{k-1})$ \hspace{7mm} $\rhd$ row indices selection via DEIM\;
        $\bm q_{\bm B_i} \gets \texttt{find\_adjacent}(\bm{q})$ \hspace{16mm} $\rhd$ Adjacent columns for $\bm B_i$\;
        $\bm p_{\bm A_i} \gets \texttt{find\_adjacent}(\bm{p})$ \hspace{16mm} $\rhd$ Adjacent rows for $\bm A_i$\;
        $\bm Z_{\bm B_i} = \left(\bm{Y}(\bm{q}, :)^{\dagger}\right)^{\top}
        \bm{Y}(\bm{q}_{B_i}, :)^\top$ \hspace{6mm} $\rhd$ Coefficient matrix for $\bm{B}_i$\;
        $\bm Z_{\bm A_i} = \bm{U}(\bm{p}_{A_i}, :) \left(\bm{U}(\bm{p}, :)\right)^\dagger$ \hspace{12mm} $\rhd$ Coefficient matrix for $\bm{A}_i$\;
        $\bm{A}_{i}^{r} = \bm{A}_{i}(\bm{p}, \bm{p}_{A_i}) \, \bm{Z}_{A_i}$ \hspace{8mm} $\rhd$ Projected $\bm{A}_i$\;
        $\bm{B}_{i}^r = \bm{Z}_{B_i} \, \bm{B}_{i}(\bm{q}_{B_i}, \bm{q})$ \hspace{8.5mm} $\rhd$ Projected $\bm{B}_i$\;
        $\bm X(:, \bm{q})$
        $\gets \texttt{GMRES\_LME}(\bm A_i, \; \bm B_i^r, \; \bm C(:, \bm{q}), \; \bm X^{k-1}(:, \bm{q}), \; \Delta t)$ \hspace{16mm} $\rhd$ Find the solution on the selected columns\;
        $\bm X(\bm{p}, :)$
        $\gets \texttt{GMRES\_LME}(\bm A_i^r, \; \bm B_i, \; \bm C(\bm{p}, :), \; \bm X^{k-1}(\bm{p}, :), \; \Delta t)$ \hspace{16mm} $\rhd$ Find the solution on the selected rows\;
        $\bm{U}^{k}$, $\bm \Sigma^{k}$, $\bm{Y}^{k}$
        $\gets \texttt{Stable\_CUR\_Algorithm}(\bm X(\bm{p}, :), \; \bm X(:, \bm{q}))$ \hspace{16mm} $\rhd$ Update the solution\;
        $\bm{U}_\Delta^{k}$, $\bm \Sigma_\Delta^{k}$, $\bm{Y}_\Delta^{k}$
        $\gets \texttt{Low\_Rank\_Approximation\_of\_$\Delta\bm X$}(\bm{U}_\Delta^{k-1}, \; \bm{Y}_\Delta^{k-1})$ \hspace{16mm} $\rhd$ Compute $\Delta\bm X$ singular values\;
        $\|\Delta\bm X\|_F = \|\bm \Sigma_\Delta^{k}\|_F = \big (\displaystyle \sum_{i=1}^{r_{\Delta}} \sigma_i^2\big )^{1/2}$\;
    }
    $\bm{U}_R^{k}$, $\bm \Sigma_R^{k}$, $\bm{Y}_R^{k}$
    $\gets \texttt{Low\_Rank\_Approximation\_of\_Residual}(\bm{U}_R^{k-1}, \; \bm{Y}_R^{k-1})$ \hspace{16mm} $\rhd$ Compute residual singular values\;
    $\|\mathbf{R}\|_F = \|\bm \Sigma_R^{k}\|_F = \big (\displaystyle \sum_{i=1}^{r_{\text{res}}} \sigma_i^2\big )^{1/2}$\;
    \If{$\|\mathbf{R}\|_F > \epsilon_{R}$}{
    $r \gets r + \Delta r$\;
    }
}
\caption{\texttt{TDB-CUR} for solving
\(
\bm{A}_0\bm X \bm{B}_0- \Delta t \sum_{i=1}^d\left( \bm{A}_i \bm{X} \bm{B}_i \right)= \bm{A}_0 \bm X^{k} \bm{B}_0  - \Delta t  \bm C
\).}
\label{alg:ITDB}
\endgroup
\end{algorithm}
\subsection{Krylov solver}
At any step of the fixed point iteration, we solve two thin matrix systems for the selected subcolumns and subrows of $\bm X$, namely for the matrices $\bm X(:,\bm q)$, Eq.~\eqref{eq:col_aux3} and $\bm X(\bm p,:)$, Eq.~\eqref{eq:row}. In this work, aiming at solving large systems, we perform these two steps by using a Krylov method. Krylov-based methods for the solution of Sylvester equations were commented in \cite{simoncini2016computational} for the resolution of LMEs, and we refer to \cite{kressner2010krylov} for a Krylov method applied to matrix equations in Kronecker product form. 

Vectorizing the unknowns in Eq.~\eqref{eq:col_aux3} and Eq.~\eqref{eq:row} leads to high-dimensional linear systems. Our objective is to carry out matrix--vector multiplications and orthogonalization steps by leveraging the Kronecker structure of the operators and the low-rank structure of the unknowns, thereby avoiding explicit vectorization of the linear systems. To explain this, let us express the Eq.~\eqref{eq:col_aux3} for solving for $\bm X(:,\bm q)$ in the form of
\[
\mathcal A_{\bm q} (\bm X^{k+1}_{\bm q}) = \mathcal B_{\bm q} ( \hat{\bm X}^{k}),
\]
where $\bm X^{k+1}_{\bm q} = \bm X(:, \bm q)$ and $\mathcal A_{\bm q}: \mathbb{R}^{n_1 r \times n_1 r} \rightarrow \mathbb{R}^{n_1 r \times n_1 r}$ and $\mathcal B_{\bm q}: \mathbb{R}^{n_1  \times n_2} \rightarrow \mathbb{R}^{n_1 \times r}$ are as follows:
\begin{subequations}
\begin{align}
\mathcal A_{\bm q} (\bm X^{k+1}_{\bm q}) &=  \bm A_0 \bm X^{k+1}_{\bm q} \bm B_{0_r} - \Delta t\big [\bm A_1\bm X^{k+1}_{\bm q}\bm B_{1_r} +  \bm A_2\bm X^{k+1}_{\bm q} \bm B_{2_r} \big ]\label{eq:kry-col1}\\
\mathcal B_{\bm q} ( \hat{\bm X}^{k}) &= \bm A_0 \bm U^{k} \bm \Sigma^k \bm Y^k(\bm q_{\bm B_0},:)^\top \bm B_0(\bm q_{\bm B_0}, \bm q)  - \Delta t  \bm C(:,\bm q), \label{eq:kry-col2}
\end{align}
\end{subequations}
where $\hat{\bm X}^{k} =  \bm U^{k} \bm \Sigma^k \bm Y^{k^\top}$ is the known solution from the previous time-step and is stored, and its subcolumns are evaluated in the factorized compressed form.  Equivalently, Eq.~\eqref{eq:col_aux3} in the vectorized form can be expressed as:
\[
 \bm A_{\bm q}  \bm x^{k+1}_{\bm q} = \bm b_{\bm q},
\]
where $\bm x^{k+1}_{\bm q} = \mbox{vec}(\bm X^{k+1}_{\bm q}) \in \mathbb{R}^{n_1r \times 1}$, and $\bm A_{\bm q} \in \mathbb{R}^{n_1 r \times n_1 r}$ and $\bm b_{\bm q} \in \mathbb{R}^{n_1r \times 1}$ are as follows: 
where
\begin{align*}
\bm A_{\bm q} &= \bm B_{0_r}^{\mathrm T} \otimes \bm A_0 - \Delta t\big [\bm B_{1_r}^{\mathrm T} \otimes \bm A_1 + \bm B_{2_r}^{\mathrm T} \otimes \bm A_2 \big ],\\
\bm b_{\bm q} &= \mbox{vec}(\mathcal B_{\bm q} (\hat{\bm X}^{k})),
\end{align*}
While one could apply a Krylov solver directly to $\bm A_{\bm q} \, \bm x^{k+1}_{\bm q} = \bm b_{\bm q}$,
we note that the action of matrix $\bm A_{\bm q}$ on $\bm x^{k+1}_{\bm q}$, i.e.,  $\bm A_{\bm q} \bm x^{k+1}_{\bm q}$, is equivalent to Eq.~\eqref{eq:kry-col1}. Hence, instead of explicitly forming $\bm A_{\bm q}$ and $\bm b_{\bm q}$, we work with the matrix Eqs.~\eqref{eq:kry-col1}--\eqref{eq:kry-col2}. In this formulation, the Krylov subspace consists of matrices of dimension $n_1 \times r$, rather than vectors of length $r n_1 \times 1$.  

Constructing the Krylov subspace requires evaluating the residual, which amounts to applying the linear operator to a matrix $\bm Q \in \mathbb{R}^{n_1 \times r}$. This is conveniently accomplished using Eqs.~\eqref{eq:kry-col1}--\eqref{eq:kry-col2}, where the residual takes the form $\bm R = \mathcal{A}_{\bm q}(\bm Q) - \mathcal{B}_{\bm q}(\hat{\bm X}^{k})$.
Moreover, the Arnoldi process requires an inner product, which in the matrix setting is given by the Frobenius inner product:
\[
\big < \bm Q_1, \bm Q_2 \big >_F = \mbox{trace}(\bm Q_1^\top \bm Q_2 ).
\]
Note that the above inner product is equal to $\mbox{vec}(\bm Q_1)^\top \mbox{vec}(\bm Q_2)$. The rest of the Krylov solver is based on the standard generalized minimal residual method (GMRES) algorithm. The GMRES iteration we used is summarised in Algorithm~\ref{alg:Arnoldi}. 

\begin{algorithm}[htbp]
\begingroup
\fontsize{9pt}{9pt}\selectfont
\SetAlgoLined
\KwIn{$\bm A_i, \; \bm B_i, \; \bm C, \; \bm{X}^{(0)}, \; \Delta t$, $m\in\mathbb{N}^*$, $\mathrm{tol}>0$}
\KwOut{$\bm{X} \in \mathbb R^{n \times r}$}
$\mathcal{A}(\bm{X}) = \bm{X} - \Delta t \sum_{i=1}^d\left( \bm{A}_i \bm{X} \bm{B}_i \right)$\;
$\mathcal{C} = \bm{X} - \Delta t C$\;
$\bm{R}^{(0)}= \mathcal{C}-\mathcal{A}(\bm{X}^{(0)})$ \hspace{20mm} $\rhd$ initial residual (in low rank format)\;
$b=\| \mathcal{C} \|_{F}$ \hspace{35.5mm} $\rhd$ the norm of the right-hand side\;
$r^{(0)}=\| \bm{R}^{(0)} \|_{F}$ \hspace{28.0mm} $\rhd$ the initial residual norm\;
$s,c = \mathrm{zeros}(m)$ \hspace{27.5mm} $\rhd$ auxiliary vectors \;
$\bm{Q}^{(0)} = \mathbf{R}^{(0)}/r^{(0)}$ \hspace{25.0mm} $\rhd$ first Krylov space element \;
$\mathcal{Q} = \left\lbrace \bm{Q}^{(0)} \right\rbrace$ \hspace{32.0mm} $\rhd$ the list of Krylov space elements \;
$H = \mathrm{zeros}(m+1,m)$ \hspace{20.0mm} $\rhd$ matrix of the least squares problem \;
$\hat{e}_1\in\mathbb{R}^{m+1}, (\hat{e}_1)_i = \delta_{1i})$ \hspace{18.0mm} $\rhd$ an auxiliary vector \;
$\beta = r^{0} \hat{e}_1$ \hspace{37.0mm} $\rhd$ least squares right-hand side\;
\For{$i=1$ \KwTo $m$}{
$H, \bm{Q}^{(k)} = \mathrm{Arnoldi}(\mathcal{A},\mathcal{Q},\mathrm{k})$ \hspace{5.0mm} $\rhd$ Arnoldi iteration \;
$H,\beta = \mathrm{Givens}(H,c,s,\beta)$ \hspace{9mm}  $\rhd$ Givens rotations to update $H,\beta$\;
\If{$\| R^{(k)} \|_F<\mathrm{tol}$}{break} 
}
$y = H(1:k,1:k)^{-1}\beta(1:k)$ \hspace{11mm}  $\rhd$ compute the coefficients\;
$\bm{X} = \bm{X}^{(0)} + \sum_{i=1}^k y_i \bm{Q}^{(i)}$ \hspace{16mm}  $\rhd$ compute the solution \;
\caption{\texttt{GMRES\_LME}}
\label{alg:Arnoldi}
\endgroup
\end{algorithm}

\subsection{Iterative CUR for steady-state LME}\label{sec:ss}

We now describe how the proposed methodology can be applied to steady-state problems. Consider the steady-state form of the generalized Sylvester equation introduced in Eq.~\eqref{eq:gse_td}, written as:
\begin{equation}\label{eq:gse_ss}
\bm A_1 \bm X \bm B_1 + \bm A_2 \bm X \bm B_2 - \bm C = 0.
\end{equation}
For small to moderate-sized LMEs, the above equation can be solved using the TDB-CUR algorithm with the following adjustment to $\mathcal A$  and $\mathcal B$ as follows:
\begin{subequations}
\begin{align}
    \mathcal A_{ss} (\bm X) &=  \bm A_1\bm X\bm B_1 +  \bm A_2\bm X\bm B_2,   \label{eq:LME_A_ss}\\
    \mathcal B_{ss}(\bm X)  &=  \bm C. \label{eq:LME_B_ss}
    \end{align}
\end{subequations}
In the above formulation, $\mathcal B_{ss}(\bm X)$ is a constant matrix; however, we retain the argument $\bm X$ to preserve the analogy between the steady-state and time-dependent cases. \Cref{alg:ITDB} can then be applied to the steady-state LME by substituting $\mathcal A \rightarrow \mathcal A_{ss}$ and $\mathcal B \rightarrow \mathcal B_{ss}$. In this interpretation, the time-step index $(k)$ plays the role of an iteration count, and the time-dependent bases naturally act as \emph{adaptive bases}, since the column and row subspaces of $\bm X$ evolve from one iteration to the next.

We make a few observations in comparing the steady-state  LMEs given by Eqs.~(\ref{eq:LME_A_ss})-(\ref{eq:LME_B_ss}) to time-dependent LMEs given by Eqs.~(\ref{eq:LME_A})-(\ref{eq:LME_B}). First we observe that as $\Delta t \rightarrow \infty$,  the LME  $\mathcal A (\bm X^{k+1}) = \mathcal B (\bm X^k)$  becomes  identical to $\mathcal A_{ss} (\bm X^{k+1}) = \mathcal B_{ss} (\bm X^k)$.

Second, for the implicit time integration of linear MDEs, the time step $\Delta t$ must remain relatively small to ensure accuracy, making the solution from the previous step a good initial guess for GMRES. As a result, GMRES can converge with a modest Krylov subspace dimension, often without requiring many restarts. In contrast, for large-scale steady-state LMEs, a poor initial guess may necessitate constructing a Krylov subspace of very large dimension, potentially exceeding the available memory. Under a fixed memory budget, this forces GMRES to perform numerous restarts, which can significantly slow convergence.

To mitigate this issue, we introduce a pseudo-time derivative and reformulate the Eq.~\eqref{eq:gse_ss} as a time-dependent problem, analogous to the Eq.~\eqref{eq:gse_td}. However, our objective here is not to resolve the transient dynamics, but rather to compute the steady-state solution, defined as the point at which the solution becomes time-invariant. The pseudo-time formulation is given by:
\begin{equation}\label{eq:gse_pt}
\dot{\bm X}(\tau) = \bm A_1 \bm X(\tau) \bm B_1 + \bm A_2 \bm X(\tau) \bm B_2 - \bm C,
\end{equation}
where $\tau$ denotes the pseudo-time variable. Since we are only interested in the steady-state solution and not in accurately resolving the transient dynamics, we choose large values for the pseudo-time step $\Delta \tau$ to accelerate convergence. Increasing $\Delta \tau$ enhances the rate of convergence but also increases the size of the Krylov subspace constructed by the GMRES solver. This occurs because larger pseudo-time steps introduce greater differences between successive iterates, requiring the solver to capture more directions in the solution space.

To balance memory limitations and convergence speed, one can gradually increase $\Delta \tau$ over iterations. This allows the initial guess to evolve smoothly toward the correct solution in the early stages. Once the iterates approach the steady state, we use a larger value of $\Delta \tau$ to drive the solution rapidly to convergence. The pseudo-time step is updated at each iteration $i$ using the following function:
\begin{equation}\label{eq:tau}
\Delta \tau = 1 + \Delta \tau_{\text{max}} \left(1 - \exp\left(-\frac{k - 1}{a}\right)\right),
\end{equation}
where $k$ is the iteration count and $\Delta \tau_{\text{max}}$ is the maximum pseudo-time step value used to accelerate the final stage of convergence and we use $\Delta \tau_{\text{max}}$ in the order of $\Delta \tau_{\text{max}}$ can be as large as $\mathcal{O}(10^3)$ and  $\mathcal{O}(10^4)$. We choose $a=25$ in all cases considered.

The difficulty of solving large-scale LMEs is further compounded by their often poor conditioning. For instance, refining the mesh typically worsens the condition number of the discrete representation of differential operators. This highlights the importance of effective preconditioning strategies for GMRES. The design of preconditioners for multi-term LMEs is an active area of research (see, e.g.,~\cite{BKR25} and references therein). Extending such approaches to TDB-CUR is certainly a worthwhile direction, but it lies beyond the scope of the present work.

\subsubsection{Rank adaptivity}
We propose a rank-adaptive strategy in which the rank is dynamically adjusted to meet a user-prescribed residual tolerance. For MDEs, rank adaptivity ensures the desired accuracy at each time step.

In addition to accuracy, rank adaptivity is also beneficial for steady-state LMEs, where it reduces the dimension of the Krylov subspace. One can start the iterative CUR with a small rank and gradually increase it until the target accuracy is achieved. This approach keeps the Krylov subspace in GMRES small, since higher-rank iterations are initialized with previously converged low-rank solutions that provide effective starting guesses. Consequently, GMRES requires only small updates, converges with a smaller Krylov basis at each stage, and maintains lower memory usage.

\subsection{Extension to generic (non-Sylvester) linear and nonlinear matrix equations}\label{sec:non_gen_syl}
In this section, we show how the presented methodology can be utilized to solve non-Sylvester multi-term LMEs. Non-Sylvester LMEs may arise, for instance, as intermediate systems when applying Newton’s method to solve nonlinear matrix equations. The range of cases that may be considered in this setting is extensive. As a representative example,  we consider LMEs in the form of 
\begin{equation}\label{eq:genMat}
 \sum_{i=1}^d\bm A_i \bm X  \bm B_i  + \mathbf{E} \circ \bm X = \bm C,   
\end{equation}
where $\circ$ denote the Hadamard (elementwise) product and $\bm{E}$ are matrices of of the same size as $\bm X$.  As we show, the Hadamard product appears in solving nonlinear matrix equations.

 Eq.~\eqref{eq:genMat} can be solved with the iterative CUR algorithm without major changes.  For example, sampling the columns can be carried as follows
\begin{equation}
\label{eq:genMat2}
\sum_{i=1}^{d} \bm{A}_i\bm{X}\bm{B}_{i}(:,\bm q) \;+\; \bm{E}(:,\bm q)\circ \bm{X}(:,\bm q) = \bm{C}(:,\bm q).
\end{equation}
Note that the Hadamard product operates element-wise and involves only the selected entries. Consequently, the term $\bm E(:,\bm q) \circ \bm X(:,\bm q)$ can be computed directly. 

\section{Demonstration cases} \label{sec:cases}

\subsection{High-order implicit time integration}
\begin{figure}
  \centering
  \begin{subfigure}[t]{0.48\textwidth}
    \centering
    \input{grid}
    \caption{2D  domain ($\Omega^{(1)}_{2D}$).}
    \label{subig:p2schem2}
  \end{subfigure}
  \begin{subfigure}[t]{0.48\textwidth}
    \centering
    \input{grid2}
    \caption{2D domain ($\Omega^{(2)}_{2D}$).}
    \label{subig:p2schem3}
  \end{subfigure}
  \begin{subfigure}[t]{0.3\textwidth}
    \centering
    \includegraphics[width=0.99\textwidth]{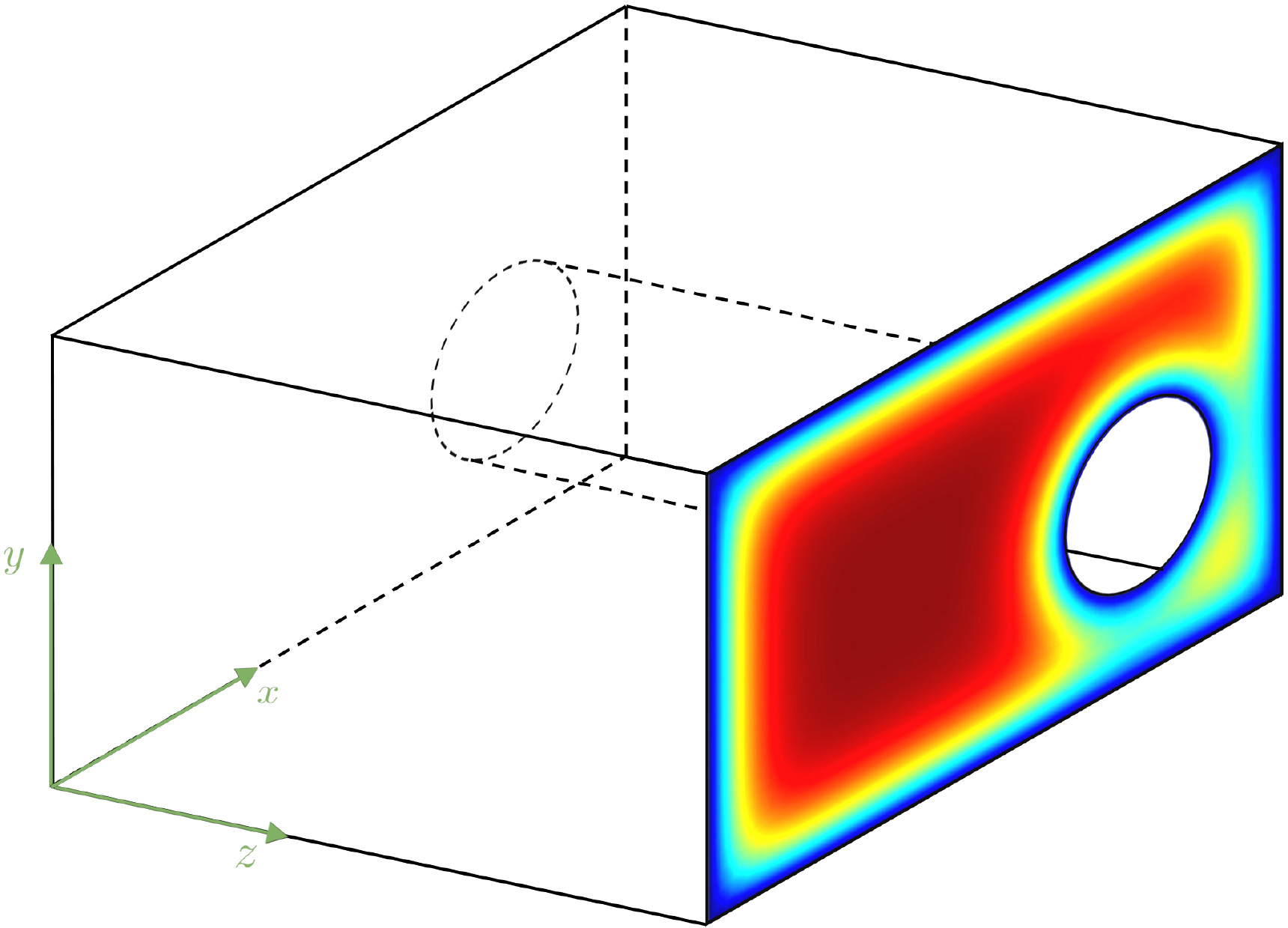}
    \caption{3D  domain ($\Omega_{3D}=\Omega^{(1)}_{2D} \times [0,1]$).}
    \label{subig:p2schem1}
  \end{subfigure}
  \caption{Computational domains used in the demonstrations: (a) two-dimensional plate with a hole denoted with $\Omega^{(1)}_{2D}$; (b) a square denoted with $\Omega^{(2)}_{2D}$; and (c) a three-dimensional domain denoted with $\Omega_{2D}$ and obtained by extruding $\Omega^{(1)}_{2D}$ in the $z$ direction.}
    \label{fig:p2schem1}
\end{figure}

As a first demonstration, we show the application of \texttt{TBD-CUR} for implicit time integration of MDEs on low-rank matrix manifolds.  To this end, we consider the three-dimensional (3D) heat equation given below:
\begin{equation}
\rho c_p \frac{\partial T}{\partial t} = k \Delta T + q, \quad (x,y,z) \in \Omega_{3D}, \quad t \in [0,0.6],
\end{equation}
 where \( \rho = 1\, \) is the density, \( c_p = 1\, \) is the specific heat capacity, \( k = 1\, \) is the thermal conductivity, and \( q = 10\} \) is a constant volumetric heat source. The temperature field is denoted by \( T(x,y,z,t) \), while \( t \in [0, 0.6] \) is the time and \( z \in [0,1] \) is the vertical coordinate.  We assume normalized (dimensionless values for all physical quantities. 
 The domain is shown schematically in Figure \ref{subig:p2schem1} and is denoted with  $\Omega_{3D}$. It is constructed via the tensor product of 
 the two-dimensional domain, shown in Figure \ref{subig:p2schem3} and denoted by \( \Omega^{(1)}_{2D} \), and  the one-dimensional domain \( z \in [0,1] \). The domain \( \Omega^{(1)}_{2D} \) is a rectangular domain \( [-1,1] \times [-0.75,0.5] \) with a circular hole centered at \( (0.5,-0.25) \) of radius \( 0.25 \) defined as:
\begin{equation}
\label{eq:omega}
\Omega^{(1)}_{2D} = \left\{ (x,y) \in [-1,1] \times [-0.75,0.5] \, \big| \, (x-0.5)^2 + (y+0.25)^2 \geq 0.25^2 \right\}.
\end{equation}

A uniform Dirichlet boundary condition is imposed on the entire spatial boundary
\(\partial\Omega_{3D} = (\partial\Omega^{(1)}_{2D} \times [0,1]) \cup (\Omega_{2D}^{(1)} \times \{0,1\})\):
\begin{equation}
T(x,y,z,t) = T_0, \qquad (x,y,z) \in \partial\Omega_{3D},\; t \in (0,0.6].
\end{equation}
The initial condition is given by:
\begin{equation}
T(x,y,z,0) = T_0, \qquad (x,y)\in \Omega^{(1)}_{2D},\; z\in[0,1],
\end{equation}
where \(T_0 = 0.5\).

We use the finite element method (FEM) for the spatial discretization and the explicit fourth-order Runge-Kutta (RK4) method for time integration to obtain the reference solution with time step $\Delta t= 1e-5$. Additionally, we solve the equation using the implicit Euler method, the second-order backward differentiation formula (BDF2), and the third-order backward differentiation formula (BDF3) using $\texttt{TDB-CUR}$. 

We now explain how the above 3D problem is matricized. We employ FEM to obtain a matrix-based formulation for this 3D problem by constructing a two-dimensional finite element basis in the \(x\)–\(y\) plane and a one-dimensional basis in the \(z\)-direction. As a result, the temperature field is stored as a matrix, where each column corresponds to a fixed \(z\)-coordinate, and each row represents the degrees of freedom associated with the \(x\)–\(y\) finite element basis. The details of the FEM model are presented in \ref{app:fem} and the resulting equation is
\begin{equation}
\label{syl_3dheat2}
\bm A_0  \frac{d \bm X}{d t} \bm B_0 = \bm A_1 \bm X \bm B_1 + \bm A_2 \bm X \bm B_2 + \bm C,
\end{equation}
 where $\mathbf{A}_i$ and $\mathbf{B}_i$ ($i=0,1,2$) represent the corresponding stiffness matrices (\(\bm K_{xy},\,\bm K_{z}\)) and the mass matrices (\(\bm M_{xy},\,\bm M_{z}\)) based on the Eq.~\eqref{eq:fem2}. 
The matrix $\mathbf{X}$ contains the unknown temperature field, and $\mathbf{C}$ accounts for heat source term.

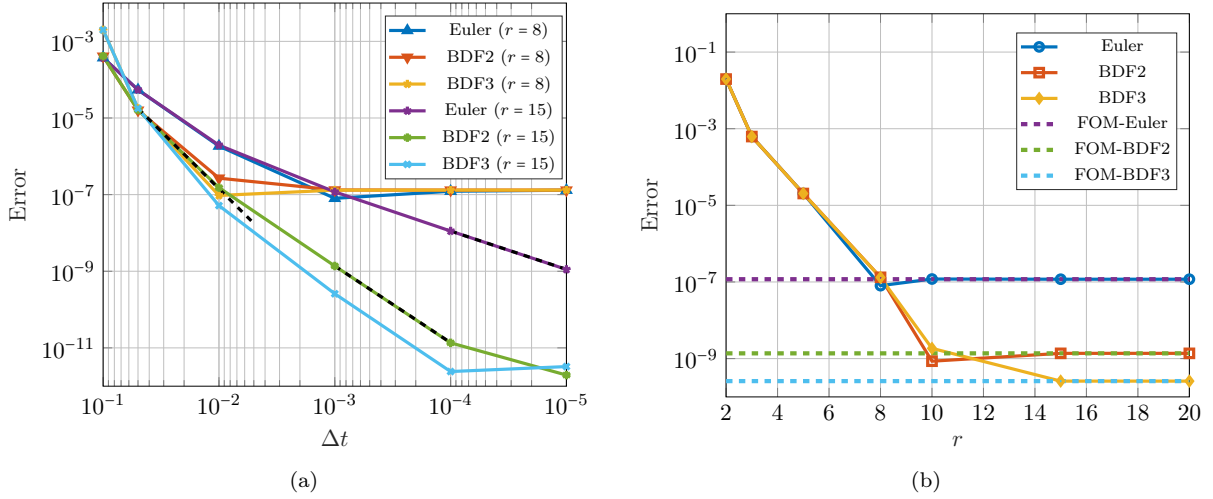
\begin{figure}[t]
  \centering
  \begin{subfigure}[t]{0.49\textwidth}
     \centering
%
%
\definecolor{mycolor1}{rgb}{0.00000,0.44700,0.74100}%
\definecolor{mycolor2}{rgb}{0.85000,0.32500,0.09800}%
\definecolor{mycolor3}{rgb}{0.92900,0.69400,0.12500}%
\definecolor{mycolor4}{rgb}{0.49400,0.18400,0.55600}%
\definecolor{mycolor5}{rgb}{0.46600,0.67400,0.18800}%
\definecolor{mycolor6}{rgb}{0.30100,0.74500,0.93300}%
\begin{tikzpicture}[scale=0.8]

\begin{axis}[%
width=3in,
height=2.5in,
at={(0.758in,0.619in)},
scale only axis,
x dir=reverse,
xmode=log,
xmin=1e-05,
xmax=0.1,
xminorticks=true,
xlabel style={font=\color{white!15!black}},
xlabel={$\Delta t$},
ymode=log,
ymin=1e-12,
ymax=0.01,
yminorticks=true,
ylabel style={font=\color{white!15!black}},
ylabel={Error},
axis background/.style={fill=white},
xmajorgrids,
xminorgrids,
ymajorgrids,
yminorgrids,
legend style={
  at={(0.8,0.54)},
  anchor=south,
  font=\scriptsize,           
  cells={align=left},           
  align=left,                   
  draw=white!15!black,          
  inner sep=1pt,                
  row sep=1pt,                  
  column sep=5pt,               
}
]
\addplot [color=mycolor1, line width=1.5pt, mark=triangle, mark options={solid, mycolor1}]
  table[row sep=crcr]{%
0.1	0.000373653127226783\\
0.05	5.43049290205466e-05\\
0.01	1.84782971576978e-06\\
0.001	8.00310106328118e-08\\
0.0001	1.22378553561619e-07\\
1e-05	1.30249485248411e-07\\
};
\addlegendentry{Euler ($r=8$)}

\addplot [color=mycolor2, line width=1.5pt, mark=triangle, mark options={solid, rotate=180, mycolor2}]
  table[row sep=crcr]{%
0.1	0.000409253824849056\\
0.05	1.59549070520143e-05\\
0.01	2.68628304831597e-07\\
0.001	1.32237758167599e-07\\
0.0001	1.31151007922461e-07\\
1e-05	1.31136651876552e-07\\
};
\addlegendentry{BDF2 ($r=8$)}

\addplot [color=mycolor3, line width=1.5pt, mark=star, mark options={solid, mycolor3}]
  table[row sep=crcr]{%
0.1	0.00198427357635249\\
0.05	1.76843402006863e-05\\
0.01	9.51073585823982e-08\\
0.001	1.30934318129595e-07\\
0.0001	1.31138196256789e-07\\
1e-05	1.31135483715952e-07\\
};
\addlegendentry{BDF3 ($r=8$)}

\addplot [color=mycolor4, line width=1.5pt, mark=star, mark options={solid, mycolor4}]
  table[row sep=crcr]{%
0.1	0.000373757171207143\\
0.05	5.44094784532977e-05\\
0.01	1.95086643353159e-06\\
0.001	1.17971477425409e-07\\
0.0001	1.12041193639798e-08\\
1e-05	1.11519787899734e-09\\
};
\addlegendentry{Euler ($r=15$)}

\addplot [color=mycolor5, line width=1.5pt, mark=asterisk, mark options={solid, mycolor5}]
  table[row sep=crcr]{%
0.1	0.000409149484547774\\
0.05	1.60594241362709e-05\\
0.01	1.52041344417297e-07\\
0.001	1.37060781791285e-09\\
0.0001	1.35223250584893e-11\\
1e-05	1.94964351464161e-12\\
};
\addlegendentry{BDF2 ($r=15$)}

\addplot [color=mycolor6, line width=1.5pt, mark=x, mark options={solid, mycolor6}]
  table[row sep=crcr]{%
0.1	0.00198437537411962\\
0.05	1.77840735213784e-05\\
0.01	5.16656324175078e-08\\
0.001	2.58431660839235e-10\\
0.0001	2.41189713341228e-12\\
1e-05	3.25028194316261e-12\\
};
\addlegendentry{BDF3 ($r=15$)}

\addplot [color=black, dashed, line width=1.5pt, forget plot]
  table[row sep=crcr]{%
0.0001	1.12041631651557e-08\\
1e-05	1.12041631651557e-09\\
};
\addplot [color=black, dashed, line width=1.5pt, forget plot]
  table[row sep=crcr]{%
0.001	1.3706218621112e-09\\
0.0001	1.3706218621112e-11\\
};
\addplot [color=black, dashed, line width=1.5pt, forget plot]
  table[row sep=crcr]{%
0.05	1.77840735537546e-05\\
0.005	1.77840735537546e-08\\
};
\node[above left, align=right, inner sep=0]
at (axis cs:0,0) {Slope = 1};
\node[above left, align=right, inner sep=0]
at (axis cs:0.001,0) {Slope = 2};
\node[above left, align=right, inner sep=0]
at (axis cs:0.05,0) {Slope = 3};
\end{axis}
\end{tikzpicture}%
    \label{fig:p2Erra}
    \caption{ }
  \end{subfigure}
  \begin{subfigure}[t]{0.49\textwidth}
    \centering
%
%
\definecolor{mycolor1}{rgb}{0.00000,0.44700,0.74100}%
\definecolor{mycolor2}{rgb}{0.85000,0.32500,0.09800}%
\definecolor{mycolor3}{rgb}{0.92900,0.69400,0.12500}%
\definecolor{mycolor4}{rgb}{0.49400,0.18400,0.55600}%
\definecolor{mycolor5}{rgb}{0.46600,0.67400,0.18800}%
\definecolor{mycolor6}{rgb}{0.30100,0.74500,0.93300}%
\begin{tikzpicture}[scale=0.8]

\begin{axis}[%
width=3in,
height=2.5in,
at={(0.758in,0.509in)},
scale only axis,
xmin=2,
xmax=20,
xlabel style={font=\color{white!15!black}},
xlabel={$r$},
ymode=log,
ymin=1e-10,
ymax=1,
yminorticks=true,
ylabel style={font=\color{white!15!black}},
ylabel={Error},
axis background/.style={fill=white},
xmajorgrids,
ymajorgrids,
yminorgrids,
legend style={
  at={(0.8,0.54)},
  anchor=south,
  font=\scriptsize,           
  cells={align=left},           
  align=left,                   
  draw=white!15!black,          
  inner sep=1pt,                
  row sep=1pt,                  
  column sep=5pt,               
}
]
\addplot [color=mycolor1, line width=1.5pt, mark=o, mark options={solid, mycolor1}]
  table[row sep=crcr]{%
2	0.019756658676176\\
3	0.000624034458479865\\
5	2.04379893080445e-05\\
8	8.00310106328118e-08\\
10	1.19359011463055e-07\\
15	1.17971477425409e-07\\
20	1.17971486137833e-07\\
};
\addlegendentry{Euler}

\addplot [color=mycolor2, line width=1.5pt, mark=square, mark options={solid, mycolor2}]
  table[row sep=crcr]{%
2	0.0197566014236972\\
3	0.000624086880679525\\
5	2.03571434800377e-05\\
8	1.32237758167599e-07\\
10	8.55060848657743e-10\\
15	1.37060781791285e-09\\
20	1.37061105598931e-09\\
};
\addlegendentry{BDF2}

\addplot [color=mycolor3, line width=1.5pt, mark=diamond, mark options={solid, mycolor3}]
  table[row sep=crcr]{%
2	0.0197566022356244\\
3	0.00062408616577311\\
5	2.03582418063728e-05\\
8	1.30934318129595e-07\\
10	1.85213171725072e-09\\
15	2.58431660839235e-10\\
20	2.58418678026046e-10\\
};
\addlegendentry{BDF3}

\addplot [color=mycolor4, dashed, line width=2.0pt]
  table[row sep=crcr]{%
2	1.17971477572231e-07\\
20	1.17971477572231e-07\\
};
\addlegendentry{FOM-Euler}

\addplot [color=mycolor5, dashed, line width=2.0pt]
  table[row sep=crcr]{%
2	1.3706218621112e-09\\
20	1.3706218621112e-09\\
};
\addlegendentry{FOM-BDF2}

\addplot [color=mycolor6, dashed, line width=2.0pt]
  table[row sep=crcr]{%
2	2.5839373478042e-10\\
20	2.5839373478042e-10\\
};
\addlegendentry{FOM-BDF3}

\end{axis}
\end{tikzpicture}%
    \label{fig:p2Errb}
    \caption{}
  \end{subfigure}
  \caption{3D FEM heat conduction on $\Omega_{3D}$:  All the errors are computed against explicit Runge Kutta fourth order with time step $\Delta t= 1e-5$. The number of grid points $n_{xy} = 462$ in the $x$ and $y$ directions combined and $n_{z} = 61$ in the $z$ direction. Final time is $T_f=0.6$. Panel (a) shows the relative error over time-step size $\Delta t$ for FOM and TDB–CUR for Euler, BDF2 and BDF3 integrators at ranks $r = 8, 15$, and panel (b) shows the relative error versus rank size $r$ for Euler, BDF2, and BDF3 integrators at time step $\Delta t = 0.001$.}
  \label{fig:p3error}
\end{figure}

\Cref{fig:p3error} shows the relative error computed with respect to the FOM reference solution obtained via explicit fourth-order Runge-Kutta (RK4). Panel (a) displays the relative error across various time-step sizes. As expected, reducing the time step improves accuracy until the error plateaus due to the low-rank approximation. Higher-order methods reach this plateau sooner (i.e., at larger values of $\Delta t$) compared to lower-order methods. Additionally, the approximation with rank $r=15$ yields lower error than with $r=8$. Panel (b) presents the error as a function of rank for a fixed time step $\Delta t=0.001$. In this panel, the error plateaus are instead governed by the time-step size rather than the rank.

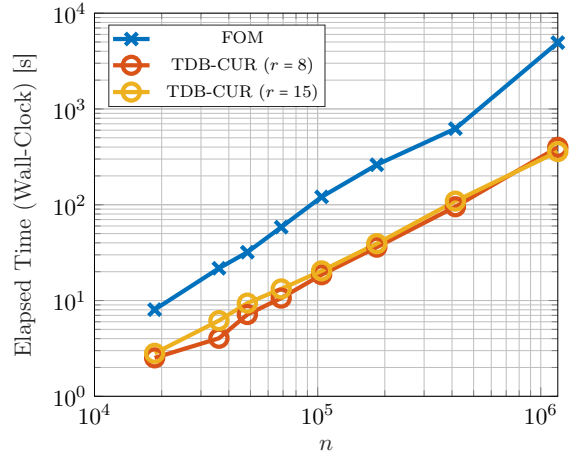
\begin{figure}[htbp]
  \centering
%
%
\definecolor{mycolor1}{rgb}{0.00000,0.44700,0.74100}%
\definecolor{mycolor2}{rgb}{0.85000,0.32500,0.09800}%
\definecolor{mycolor3}{rgb}{0.92900,0.69400,0.12500}%
\begin{tikzpicture}[scale=0.8]

\begin{axis}[%
width=3in,
height=2.5in,
at={(0.758in,0.534in)},
scale only axis,
xmode=log,
xmin=10000,
xmax=1192002,
xminorticks=true,
xlabel style={font=\color{white!15!black}},
xlabel={$n$},
ymode=log,
ymin=1,
ymax=10000,
yminorticks=true,
ylabel style={font=\color{white!15!black}},
ylabel={Elapsed Time (Wall-Clock) [s]},
axis background/.style={fill=white},
xmajorgrids,
xminorgrids,
ymajorgrids,
yminorgrids,
legend style={
  at={(0.03,0.97)},
  anchor=north west,
  font=\scriptsize,           
  cells={align=left},           
  align=left,                   
  draw=white!15!black,          
  inner sep=1pt,                
  row sep=1pt,                  
  column sep=5pt,               
}
]
\addplot [color=mycolor1, line width=2.0pt, mark size=4.0pt, mark=x, mark options={solid, mycolor1}]
  table[row sep=crcr]{%
18605	8.032242\\
36112	21.706277\\
48434	31.931351\\
68686	58.322953\\
104188	120.74793\\
184159	261.459809\\
414434	622.4928\\
1192002	4924.487553\\
};
\addlegendentry{FOM}

\addplot [color=mycolor2, line width=2.0pt, mark size=4.0pt, mark=o, mark options={solid, mycolor2}]
  table[row sep=crcr]{%
18605	2.52071\\
36112	4.035394\\
48434	7.139581\\
68686	10.603162\\
104188	18.585704\\
184159	36.050064\\
414434	95.310103\\
1192002	397.814077\\
};
\addlegendentry{TDB-CUR ($r=8$)}

\addplot [color=mycolor3, line width=2.0pt, mark size=4.0pt, mark=o, mark options={solid, mycolor3}]
  table[row sep=crcr]{%
18605	2.812895\\
36112	6.169053\\
48434	9.39499\\
68686	13.234028\\
104188	20.397211\\
184159	39.264855\\
414434	109.064626\\
1192002	358.854621\\
};
\addlegendentry{TDB-CUR ($r=15$)}

\end{axis}
\end{tikzpicture}%
  \caption{3D FEM heat conduction on $\Omega_{3D}$: the figure shows elapsed time for FOM and TDB-CUR methods at rank $r=8, 15$. Elapsed time, measured using MATLAB’s tic/toc, reflects total wall-clock time. The horizontal axis $n$ is the total number of grid points, i.e., $n=n_{xy} n_z$.}
  \label{fig:p2cputime}
\end{figure}

\Cref{fig:p2cputime} shows the wall-clock time (elapsed time) required to solve the 3D heat conduction problem using FOM and TDB-CUR with ranks \( r = 8 \) and \( r = 15 \). The figure demonstrates significant reductions in wall-clock time achieved by the reduced-order models as the system size $n$ increases, 
where $n = n_{xy} n_{z}$ denotes the total number of grid points. Therefore, increasing $n$ corresponds to refining the spatial resolution in both directions. For example,  the TDB-CUR method with \( r = 15 \) solves a system with over \( 10^6 \) unknowns in approximately the same time the FOM requires to solve a system of size around \( 2 \times 10^4 \), $50$ times faster, highlighting the efficiency and scalability of the proposed method for large-scale problems. 

\begin{figure}
  \centering

  \begin{subfigure}[t]{0.32\textwidth}
    \centering
    \includegraphics[width=0.99\textwidth]{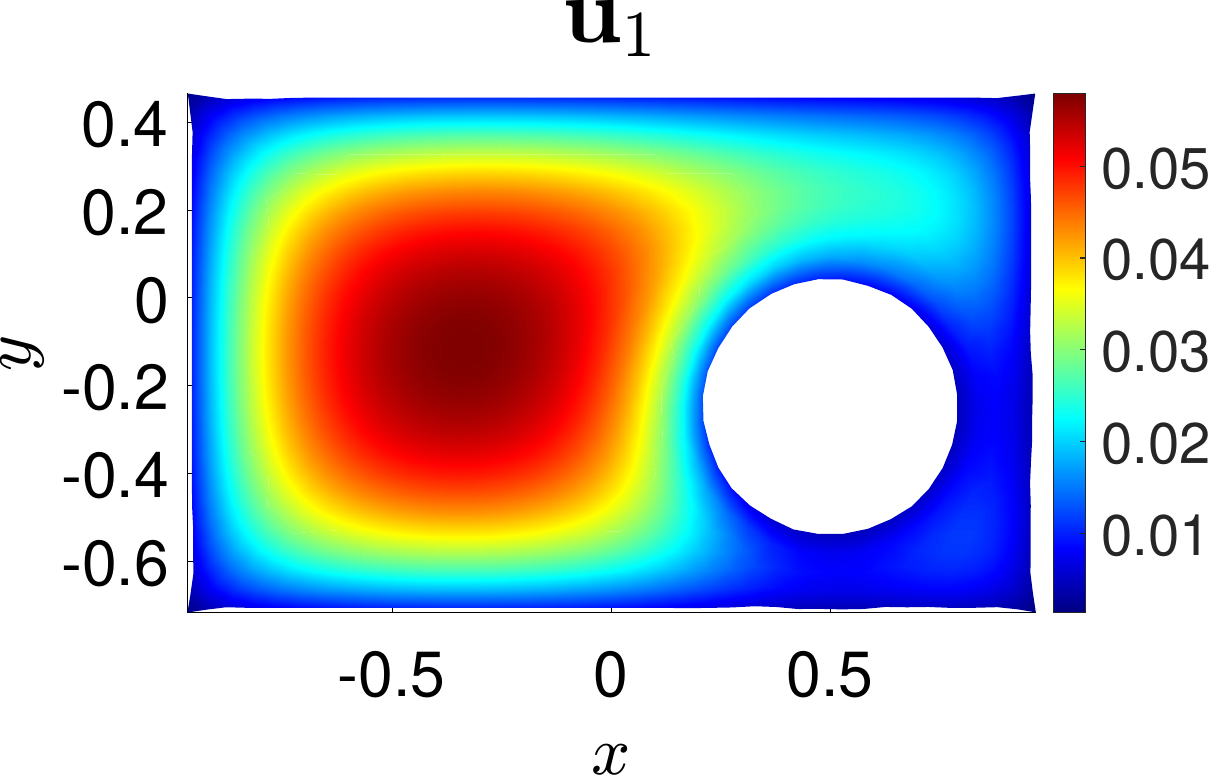}
  \end{subfigure}
  \hfill
  \hspace*{0.1cm}
  \begin{subfigure}[t]{0.32\textwidth}
    \centering
    \includegraphics[width=0.99\textwidth]{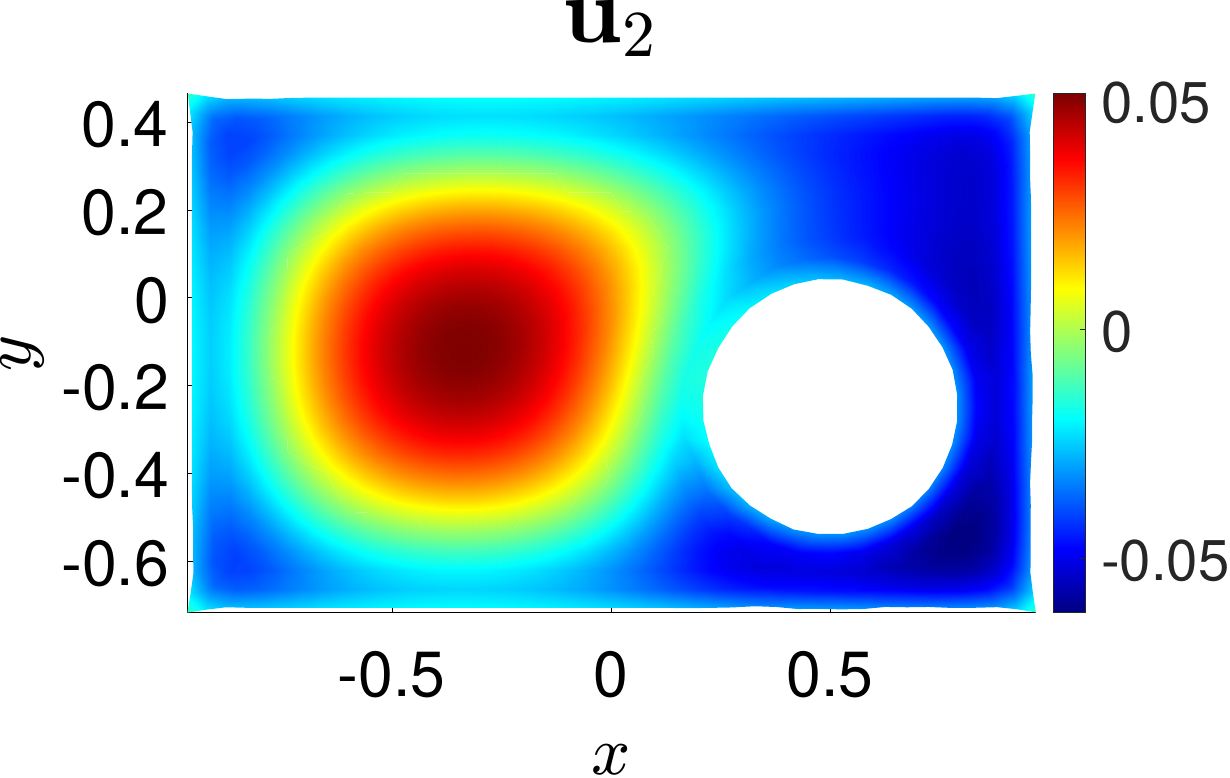}
  \end{subfigure}
  \hfill
  \hspace*{0.1cm}
  \begin{subfigure}[t]{0.32\textwidth}
    \centering
    \includegraphics[width=0.99\textwidth]{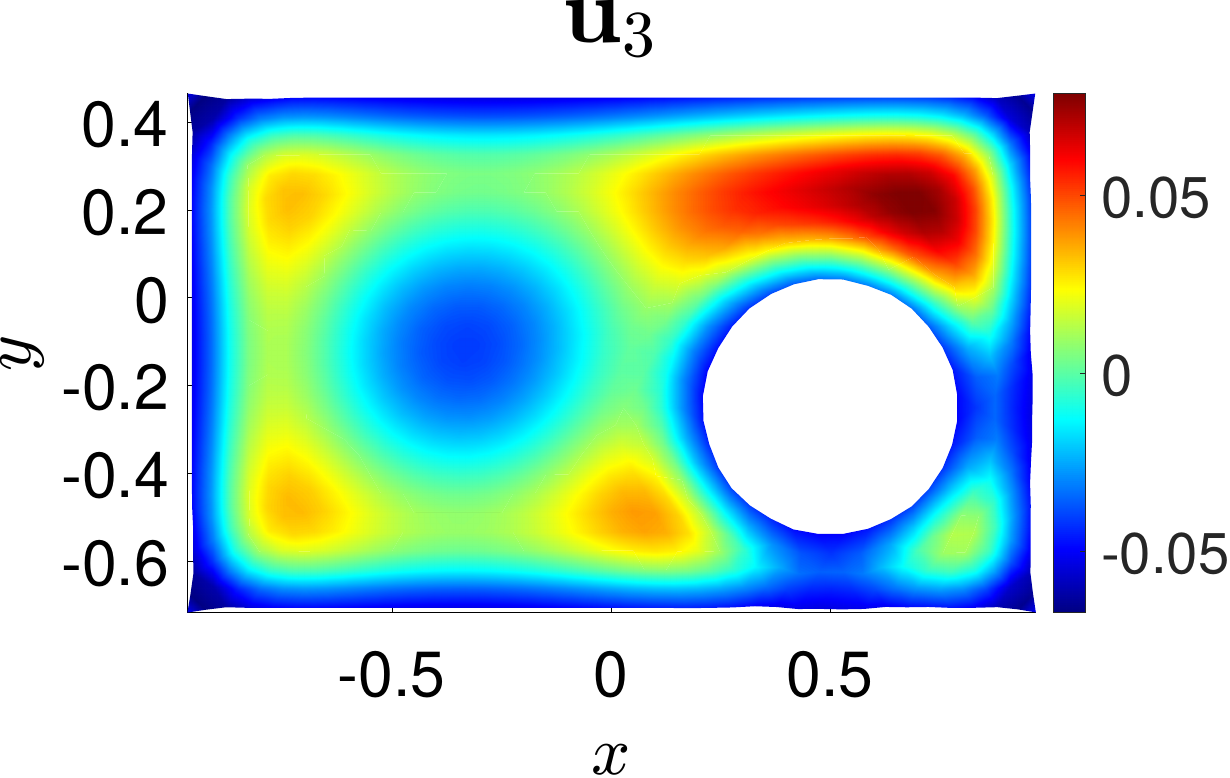}
  \end{subfigure}
  
  \vskip\baselineskip

  \makebox[\textwidth][c]{
  \hspace*{-1.0cm}
  \begin{subfigure}[t]{0.37\textwidth}
    \centering
%
%
\definecolor{mycolor1}{rgb}{0.00000,0.44700,0.74100}%
\begin{tikzpicture}

\begin{axis}[%
width=1.7in,
height=0.8631in,
at={(0.758in,1.324in)},
scale only axis,
separate axis lines,
every outer x axis line/.append style={black},
every x tick label/.append style={font=\color{black}},
every x tick/.append style={black},
xmin=0.01,
xmax=0.99,
xlabel={$z$},
every outer y axis line/.append style={black},
every y tick label/.append style={font=\color{black}},
every y tick/.append style={black},
ymin=-0.14227114888476,
ymax=0.265421308205912,
axis background/.style={fill=white},
title style={font=\bfseries},
title={$\mathbf{y}_1$},
xmajorgrids,
ymajorgrids
]
\addplot [color=mycolor1, , line width=1.2pt, forget plot]
  table[row sep=crcr]{%
0.01	0.209506475630352\\
0.02	0.265421308205912\\
0.03	0.240681457361089\\
0.04	0.178641189750716\\
0.05	0.104659902174552\\
0.06	0.0330450992456353\\
0.07	-0.0287034547108655\\
0.08	-0.0772573809807823\\
0.09	-0.111802609565875\\
0.1	-0.132990274420588\\
0.11	-0.14227114888476\\
0.12	-0.141475162335017\\
0.13	-0.132549029728658\\
0.14	-0.11739651505675\\
0.15	-0.0977854727868651\\
0.16	-0.075298246333513\\
0.17	-0.0513099995697644\\
0.18	-0.0269847665797765\\
0.19	-0.00328243538852394\\
0.2	0.0190278387542473\\
0.21	0.0393507688014643\\
0.22	0.057246571182625\\
0.23	0.0724123584655958\\
0.24	0.0846640780706177\\
0.25	0.093919469233286\\
0.26	0.100182297846537\\
0.27	0.103527987347274\\
0.28	0.104090671098746\\
0.29	0.102051632446761\\
0.3	0.0976290622396388\\
0.31	0.0910690428096815\\
0.32	0.0826376571568554\\
0.33	0.0726141188424542\\
0.34	0.0612848194573221\\
0.35	0.0489381947866327\\
0.36	0.0358603167896818\\
0.37	0.0223311254544386\\
0.38	0.00862122193657586\\
0.39	-0.00501084820491035\\
0.4	-0.0183208856330205\\
0.41	-0.0310809587603972\\
0.42	-0.0430810066928024\\
0.43	-0.054130167996926\\
0.44	-0.0640578740255199\\
0.45	-0.0727147402664914\\
0.46	-0.0799732843052163\\
0.47	-0.0857284944567412\\
0.48	-0.0898982688987641\\
0.49	-0.0924237411738278\\
0.5	-0.0932695041780947\\
0.51	-0.0924237411738076\\
0.52	-0.0898982688987969\\
0.53	-0.0857284944567054\\
0.54	-0.0799732843052624\\
0.55	-0.0727147402664713\\
0.56	-0.0640578740255007\\
0.57	-0.0541301679968671\\
0.58	-0.0430810066927304\\
0.59	-0.0310809587603602\\
0.6	-0.0183208856329711\\
0.61	-0.00501084820485664\\
0.62	0.0086212219366515\\
0.63	0.0223311254544549\\
0.64	0.035860316789708\\
0.65	0.0489381947866807\\
0.66	0.061284819457343\\
0.67	0.0726141188425046\\
0.68	0.0826376571569046\\
0.69	0.0910690428097035\\
0.7	0.0976290622396817\\
0.71	0.102051632446814\\
0.72	0.104090671098796\\
0.73	0.103527987347336\\
0.74	0.100182297846571\\
0.75	0.0939194692333394\\
0.76	0.0846640780706567\\
0.77	0.0724123584656442\\
0.78	0.0572465711826482\\
0.79	0.0393507688014946\\
0.8	0.019027838754296\\
0.81	-0.00328243538849497\\
0.82	-0.0269847665797704\\
0.83	-0.0513099995697369\\
0.84	-0.0752982463334999\\
0.85	-0.0977854727868638\\
0.86	-0.117396515056747\\
0.87	-0.132549029728656\\
0.88	-0.141475162335019\\
0.89	-0.142271148884757\\
0.9	-0.132990274420609\\
0.91	-0.111802609565891\\
0.92	-0.0772573809807799\\
0.93	-0.0287034547108524\\
0.94	0.0330450992456316\\
0.95	0.104659902174556\\
0.96	0.178641189750715\\
0.97	0.24068145736109\\
0.98	0.265421308205906\\
0.99	0.20950647563033\\
};
\end{axis}
\end{tikzpicture}%
  \end{subfigure}
  \hfill
  \hspace*{-1.0cm}
  \begin{subfigure}[t]{0.37\textwidth}
    \centering
%
%
\definecolor{mycolor1}{rgb}{0.00000,0.44700,0.74100}%
\begin{tikzpicture}

\begin{axis}[%
width=1.7in,
height=0.8631in,
at={(0.758in,1.174in)},
scale only axis,
separate axis lines,
every outer x axis line/.append style={black},
every x tick label/.append style={font=\color{black}},
every x tick/.append style={black},
xmin=0.01,
xmax=0.99,
xlabel={$z$},
every outer y axis line/.append style={black},
every y tick label/.append style={font=\color{black}},
every y tick/.append style={black},
ymin=-0.172672110047465,
ymax=0.300078560371751,
axis background/.style={fill=white},
title style={font=\bfseries},
title={$\mathbf{y}_2$},
xmajorgrids,
ymajorgrids
]
\addplot [color=mycolor1, , line width=1.2pt, forget plot]
  table[row sep=crcr]{%
0.01	0.300078560371751\\
0.02	0.256564469670365\\
0.03	0.113658446830629\\
0.04	-0.0239241593139173\\
0.05	-0.119137257685499\\
0.06	-0.166264490454328\\
0.07	-0.17267211004743\\
0.08	-0.149857962601594\\
0.09	-0.109281842049714\\
0.1	-0.0606343142236933\\
0.11	-0.0113200855321827\\
0.12	0.0334838100430147\\
0.13	0.0705266300099369\\
0.14	0.0981075042270004\\
0.15	0.115705363467068\\
0.16	0.123650885321223\\
0.17	0.122852507116272\\
0.18	0.114576835269445\\
0.19	0.100278428024134\\
0.2	0.0814719270191782\\
0.21	0.059639216057554\\
0.22	0.0361647987536632\\
0.23	0.0122934318463322\\
0.24	-0.0108950212988114\\
0.25	-0.0324976649816115\\
0.26	-0.0517909281338303\\
0.27	-0.0682258778465635\\
0.28	-0.0814185935073568\\
0.29	-0.0911372373594371\\
0.3	-0.0972870776854553\\
0.31	-0.0998944163715406\\
0.32	-0.0990901306967187\\
0.33	-0.095093348561939\\
0.34	-0.0881956274483729\\
0.35	-0.0787458920140573\\
0.36	-0.0671362968025127\\
0.37	-0.053789113487255\\
0.38	-0.0391446920289653\\
0.39	-0.0236505084656999\\
0.4	-0.00775128593106031\\
0.41	0.00811984239999964\\
0.42	0.0235491712472791\\
0.43	0.0381493113653691\\
0.44	0.0515651799497101\\
0.45	0.0634791451536743\\
0.46	0.0736153741492997\\
0.47	0.0817434300761193\\
0.48	0.0876811578734878\\
0.49	0.0912968927592499\\
0.5	0.0925110181871048\\
0.51	0.0912968927591703\\
0.52	0.0876811578739244\\
0.53	0.0817434300763933\\
0.54	0.0736153741493469\\
0.55	0.0634791451536894\\
0.56	0.0515651799499753\\
0.57	0.0381493113655161\\
0.58	0.0235491712471039\\
0.59	0.00811984240004789\\
0.6	-0.00775128593132161\\
0.61	-0.023650508465534\\
0.62	-0.0391446920290839\\
0.63	-0.0537891134870596\\
0.64	-0.0671362968024369\\
0.65	-0.0787458920141228\\
0.66	-0.0881956274485893\\
0.67	-0.0950933485619565\\
0.68	-0.0990901306966791\\
0.69	-0.0998944163717238\\
0.7	-0.0972870776853066\\
0.71	-0.0911372373593556\\
0.72	-0.0814185935073971\\
0.73	-0.0682258778469375\\
0.74	-0.0517909281337942\\
0.75	-0.032497664981428\\
0.76	-0.0108950212986834\\
0.77	0.0122934318463684\\
0.78	0.0361647987537937\\
0.79	0.0596392160575301\\
0.8	0.0814719270188724\\
0.81	0.100278428023818\\
0.82	0.114576835269713\\
0.83	0.122852507116047\\
0.84	0.123650885321156\\
0.85	0.115705363467134\\
0.86	0.098107504226975\\
0.87	0.0705266300100921\\
0.88	0.0334838100430797\\
0.89	-0.011320085531926\\
0.9	-0.0606343142237592\\
0.91	-0.109281842049436\\
0.92	-0.149857962601653\\
0.93	-0.172672110047465\\
0.94	-0.166264490454305\\
0.95	-0.119137257685505\\
0.96	-0.0239241593139987\\
0.97	0.113658446830578\\
0.98	0.25656446967054\\
0.99	0.300078560371397\\
};
\end{axis}
\end{tikzpicture}%
  \end{subfigure}
  \hfill
  \hspace*{-1.0cm}
  \begin{subfigure}[t]{0.37\textwidth}
    \centering
%
%
\definecolor{mycolor1}{rgb}{0.00000,0.44700,0.74100}%
\begin{tikzpicture}

\begin{axis}[%
width=1.7in,
height=0.8631in,
at={(0.758in,1.026in)},
scale only axis,
separate axis lines,
every outer x axis line/.append style={black},
every x tick label/.append style={font=\color{black}},
every x tick/.append style={black},
xmin=0.01,
xmax=0.99,
xlabel={$z$},
every outer y axis line/.append style={black},
every y tick label/.append style={font=\color{black}},
every y tick/.append style={black},
ymin=-0.34055927067545,
ymax=0.196037792311864,
axis background/.style={fill=white},
title style={font=\bfseries},
title={$\mathbf{y}_3$},
xmajorgrids,
ymajorgrids
]
\addplot [color=mycolor1, line width=1.2pt, forget plot]
  table[row sep=crcr]{%
0.01	-0.34055927067545\\
0.02	-0.127165274984449\\
0.03	0.0954076306091527\\
0.04	0.196037792311864\\
0.05	0.186923330347678\\
0.06	0.115195665699184\\
0.07	0.0240230686815636\\
0.08	-0.0577364832949325\\
0.09	-0.115173583459233\\
0.1	-0.143767720290914\\
0.11	-0.145574753299275\\
0.12	-0.126140895457113\\
0.13	-0.0923206975765042\\
0.14	-0.0508856366589606\\
0.15	-0.00773448183720166\\
0.16	0.0324752488900717\\
0.17	0.0664208686505713\\
0.18	0.0920416039528298\\
0.19	0.108375213349669\\
0.2	0.115354954540055\\
0.21	0.11359921488573\\
0.22	0.104213100326796\\
0.23	0.0886124210002065\\
0.24	0.0683746874467222\\
0.25	0.0451180491165907\\
0.26	0.020406900500957\\
0.27	-0.00431833385069374\\
0.28	-0.0277902879753714\\
0.29	-0.0489458210163251\\
0.3	-0.0669423676028562\\
0.31	-0.0811623391071503\\
0.32	-0.0912081343556997\\
0.33	-0.0968899863978283\\
0.34	-0.0982085443757196\\
0.35	-0.0953337828275858\\
0.36	-0.0885815532723947\\
0.37	-0.0783888485436664\\
0.38	-0.0652886388477097\\
0.39	-0.049884959076992\\
0.4	-0.0328287765589113\\
0.41	-0.0147950442108967\\
0.42	0.00353875689571584\\
0.43	0.0215123627343142\\
0.44	0.0385006003817652\\
0.45	0.0539292815885715\\
0.46	0.0672890439516063\\
0.47	0.0781470952107487\\
0.48	0.0861568354566937\\
0.49	0.0910653449140678\\
0.5	0.0927187325644767\\
0.51	0.0910653449136085\\
0.52	0.0861568354571539\\
0.53	0.0781470952103043\\
0.54	0.067289043952564\\
0.55	0.0539292815880078\\
0.56	0.0385006003831028\\
0.57	0.0215123627351346\\
0.58	0.00353875689578832\\
0.59	-0.014795044211224\\
0.6	-0.0328287765594565\\
0.61	-0.0498849590775836\\
0.62	-0.0652886388483294\\
0.63	-0.0783888485436217\\
0.64	-0.0885815532711325\\
0.65	-0.0953337828257716\\
0.66	-0.0982085443761274\\
0.67	-0.0968899863953388\\
0.68	-0.0912081343542496\\
0.69	-0.0811623391060847\\
0.7	-0.0669423676011267\\
0.71	-0.0489458210151486\\
0.72	-0.0277902879747035\\
0.73	-0.00431833384935649\\
0.74	0.0204069005016622\\
0.75	0.0451180491170285\\
0.76	0.0683746874468626\\
0.77	0.0886124210013214\\
0.78	0.104213100326996\\
0.79	0.113599214885535\\
0.8	0.115354954540193\\
0.81	0.108375213348623\\
0.82	0.092041603953232\\
0.83	0.0664208686491764\\
0.84	0.0324752488899514\\
0.85	-0.0077344818383403\\
0.86	-0.0508856366596597\\
0.87	-0.0923206975779004\\
0.88	-0.126140895457068\\
0.89	-0.145574753299588\\
0.9	-0.143767720291308\\
0.91	-0.115173583458807\\
0.92	-0.0577364832949365\\
0.93	0.0240230686807727\\
0.94	0.115195665698384\\
0.95	0.186923330347195\\
0.96	0.196037792311454\\
0.97	0.095407630608782\\
0.98	-0.127165274983908\\
0.99	-0.340559270670039\\
};
\end{axis}
\end{tikzpicture}%
  \end{subfigure}
  }
  
  \caption{First three modes from low-rank decomposition of the solution. The top row shows the first three left singular vectors \( \mathbf{u}_1, \mathbf{u}_2, \mathbf{u}_3 \), representing spatial modes in the \( x \)-\( y \) plane. The bottom row shows the corresponding right singular vectors \( \mathbf{y}_1, \mathbf{y}_2, \mathbf{y}_3 \), which capture the spatial variation along the \( z \)-direction. All the plots are for $t=0.6$s. 
}
  \label{fig:p2modes}
\end{figure}

\Cref{fig:p2modes} shows the first three dominant modes obtained from the low-rank decomposition of the solution. The top row presents the spatial modes (left singular vectors) in the \( x \)-\( y \) plane as contour plots, while the bottom row displays the corresponding spatial modes (right singular vectors) along the \( z \) direction. Note that these singular vectors are derived by expressing the CUR decomposition in the SVD form generated by \Cref{alg:SCUR}. These dominant modes illustrate how the temperature field is distributed throughout the domain.

\subsection{Non-Sylvester matrix equation}
In this demonstration, we show the application of the presented  methodology for solving nonlinear matrix equations using Newton's method. In particular, we show that the linearzing the nonlinear equation results in non-Sylvester LME. 
To this end, we solve the steady-state three-dimensional (3D) heat equation on the same geometry as in the previous section, but instead of applying a volumetric heat source throughout the domain, we impose a radiative heat flux on the surface of a circular region, leading to a nonlinear PDE. The governing equation becomes:

Let the domain be $D =\Omega\times [0,1]\subset \mathbb{R}^3$ and let $\Gamma_1 = \left\lbrace (x,y,z)\in\mathbb{R}^3 : \ (x-0.5)^2 + (y+0.25)^2 = 0.25 \right\rbrace$. The governing equation is:
\begin{equation} \label{eq:steady_heat_radiation}
\begin{cases}
\displaystyle -\Delta T = 0, &  \mathrm{in} \ \Omega\times [0,1], \\[10pt]
\displaystyle - \frac{\partial T}{\partial n} = \epsilon \sigma (T^4 - T_{\infty}^4), & \text{on } \Gamma_1, \\
T = 313.15, &  \mathrm{on} \ \partial D \setminus \Gamma_1. 
\end{cases}
\end{equation}
where \(\Omega\) is defined in Eq.~\eqref{eq:omega}. The radiation boundary condition is applied only within the circular region. Here, \( \epsilon = 0.9 \) is the surface emissivity, \( \sigma = 5.67 \times 10^{-8}\,\mathrm{W/(m^2 \cdot K^4)} \) is the Stefan–Boltzmann constant, and \( T_{\infty} = 273.15 \) is the ambient temperature. On the rectangle edge, a Dirichlet boundary condition is imposed.
The details of the FEM model are presented in \ref{app:fem_rad}. The resulting equations to be solved are as follows: 
\begin{equation}\label{eq:newton11}
\bm{A}_1 \delta \bm{X} \bm{B}_1 + \bm{A}_2 \delta \bm{X} \bm{B}_2 
+ 4 \delta \epsilon \, \bm{G} \circ \left((\bm{X}^k)^3 \circ \delta \bm{X}\right)
= -\mathcal{R}(\bm{X}^k),
\end{equation}
where \(\mathcal{R}(\bm{X})\) is the residual 
\[
\mathcal{R}(\bm{X}) = \bm{A}_{1} \bm{X} \bm{B}_1 + \bm{A}_{2} \bm{X} \bm{B}_2 
- \epsilon \sigma \left[ \bm{G} \circ \left( \bm{X}^4 - \bm{X}_\infty^4 \right)\right],
\]
and \(\delta \bm{X}\) is the unknown.

Eq.~\eqref{eq:newton11} is a steady-state LME, and we recast these equations in the form $\mathcal A_{ss} (\delta \bm X^{k+1}) = \mathcal B_{ss} (\bm X^k)$ as described in \Cref{sec:ss}, where $k$ is the iteration count. No pseudo-time term is added here.

\begin{figure}
  \centering
  \begin{subfigure}[t]{0.48\textwidth}
    \centering
    \input{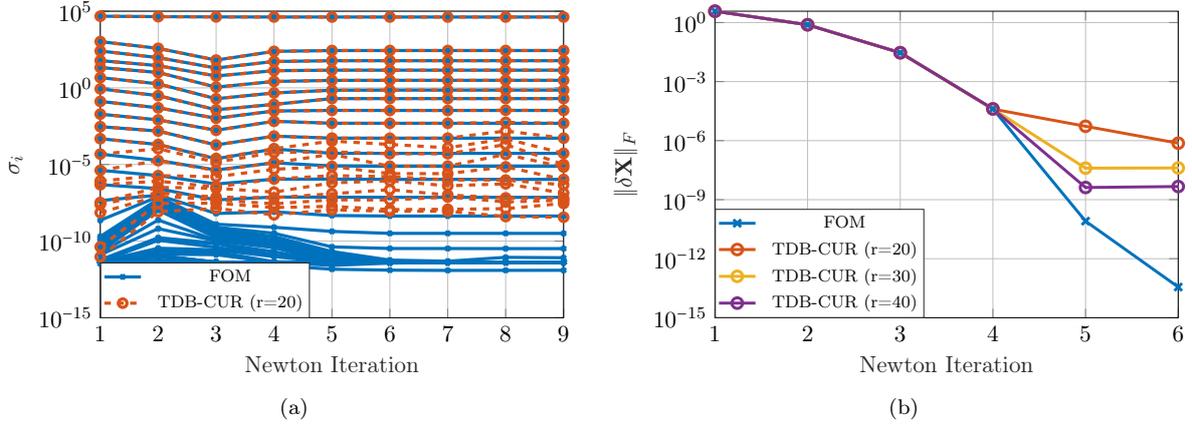}
      \label{fig:pnonlinear_a}
      \caption{}
  \end{subfigure}
  \begin{subfigure}[t]{0.48\textwidth}
    \centering
    \definecolor{mycolor1}{rgb}{0.00000,0.44700,0.74100}%

\definecolor{mycolor4}{rgb}{0.85000,0.32500,0.09800}%
\definecolor{mycolor5}{rgb}{0.92900,0.69400,0.12500}%
\definecolor{mycolor6}{rgb}{0.49400,0.18400,0.55600}%

\begin{tikzpicture}[scale=0.8]

\begin{axis}[%
width=3in,
height=2in,
at={(0.758in,0.481in)},
scale only axis,
xmin=1,
xmax=6,
xlabel style={font=\color{white!15!black}},
xlabel={Newton Iteration},
ymode=log,
ymin=1e-15,
ymax=3.7702159390128,
yminorticks=true,
ylabel style={font=\color{white!15!black}},
ylabel={$\| \delta \bm{X} \|_F$ },
axis background/.style={fill=white},
xmajorgrids,
ymajorgrids,
yminorgrids,
legend style={
  at={(0,0.355)},
  anchor=north west,
  font=\scriptsize,           
  cells={align=left},           
  align=left,                   
  draw=white!15!black,          
  inner sep=1pt,                
  row sep=1pt,                  
  column sep=5pt,               
}
]
\addplot [color=mycolor1, line width=1.4pt, mark size=2.5pt, mark=x, mark options={solid, mycolor1}]
  table[row sep=crcr]{%
1	3.7702159390128\\
2	0.770682554931008\\
3	0.0292987067346094\\
4	4.11069492337567e-05\\
5	8.09870095262769e-11\\
6	3.6357856358653e-14\\
};
\addlegendentry{FOM}



\addplot [color=mycolor4, line width=1.4pt, mark size=2.5pt, mark=o, mark options={solid, mycolor4}]
  table[row sep=crcr]{%
1	3.77021593935726\\
2	0.770682728930834\\
3	0.0292984262625188\\
4	4.00557224672158e-05\\
5	5.36292886576562e-06\\
6	7.40495826768337e-07\\
7	9.29144364073052e-07\\
};
\addlegendentry{TDB-CUR (r=20)}

\addplot [color=mycolor5, line width=1.4pt, mark size=2.5pt, mark=o, mark options={solid, mycolor5}]
  table[row sep=crcr]{%
1	3.77021593897244\\
2	0.770682555021666\\
3	0.0292986817064579\\
4	4.11349170085574e-05\\
5	3.98787724552069e-08\\
6	4.09238686945669e-08\\
7	4.9680851278706e-08\\
};
\addlegendentry{TDB-CUR (r=30)}

\addplot [color=mycolor6, line width=1.4pt, mark size=2.5pt, mark=o, mark options={solid, mycolor6}]
  table[row sep=crcr]{%
1	3.77021593929585\\
2	0.770682554676521\\
3	0.0292987079376278\\
4	4.11046548823556e-05\\
5	4.14824612340931e-09\\
6	4.65309159288517e-09\\
7	9.50170406829259e-10\\
};
\addlegendentry{TDB-CUR (r=40)}

\end{axis}
\end{tikzpicture}%
     \label{fig:pnonlinear_b}
     \caption{}
  \end{subfigure}
  \caption{Steady-state heat  conduction-radiation on $\Omega^{(1)}_{2D}$ with  $n = 18620$. Panel (a) shows the singular values of the solution over the number of Newton iterations for TDB-CUR with  ($r=20$) and all the singular values of FOM. The right panel shows the Frobenius norm of the Newton discrepancy matrix of the TDB-CUR solver for multiple values of rank and FOM.}
  \label{fig:pnonlinear}
\end{figure}

 \Cref{fig:pnonlinear_a} reports the singular values versus Newton iterations for both the FOM and TDB-CUR (\(r=20\)).  The FOM is solved by vectorizing the LME in Eq.~(\ref{eq:newton11});  after each Newton iteration, the solution is matricized and its SVD is computed. In contrast, TDB-CUR directly solves Eq.~(\ref{eq:newton11}) in low-rank form within each Newton iteration.

We observe that the TDB-CUR singular values quickly converge to $r$ leading singular values of FOM. \Cref{fig:pnonlinear_b} shows the Frobenius norm of the Newton correction, \(\lVert \delta \bm{X} \rVert_F\), versus iteration for the FOM and for TDB-CUR at several ranks. The solver terminates when \(\lVert \delta \bm{X} \rVert_F\) falls below a prescribed tolerance.  As expected, increasing the rank yields smaller values of \(\lVert \delta \bm{X} \rVert_F\) and hence higher accuracy, whereas lower ranks exhibit stagnation at larger values due to the low-rank approximation error.

\begin{figure}
  \centering
  
  \begin{subfigure}[t]{0.23\textwidth}
    \centering
%
%
\definecolor{mycolor1}{rgb}{0.00000,0.44700,0.74100}%
\definecolor{mycolor2}{rgb}{0.85000,0.32500,0.09800}%
\begin{tikzpicture}[scale=0.7]

\begin{axis}[%
width=1.5in,
height=1in,
at={(0.758in,0.481in)},
scale only axis,
unbounded coords=jump,
xmin=1,
xmax=6,
xlabel style={font=\color{white!15!black}},
xlabel={Newton Iteration},
ymin=41500,
ymax=42500,
yminorticks=true,
ylabel style={font=\color{white!15!black}, align=center},
ylabel={$\sigma_1$},
axis background/.style={fill=white},
xmajorgrids,
ymajorgrids,
yminorgrids,
legend style={
  at={(0.05,0.95)},
  anchor=north west,
  font=\scriptsize,           
  cells={align=left},           
  align=left,                   
  draw=white!15!black,          
  inner sep=1pt,                
  row sep=1pt,                  
  column sep=5pt,               
}
]
\addplot [color=mycolor1, line width=1.5pt, mark size=1.5pt, mark=x, mark options={solid, mycolor1}]
  table[row sep=crcr]{%
1	41925.030017804\\
2	41764.1965764378\\
3	41758.1045061143\\
4	41758.0959502816\\
5	41758.0959502648\\
6	41758.0959502648\\
};
\addlegendentry{FOM}
\addplot [color=mycolor2, dashed, line width=1.5pt, mark size=2.0pt, mark=o, mark options={solid, mycolor2}]
  table[row sep=crcr]{%
1	41925.0300183921\\
2	41764.1964126367\\
3	41758.1046745962\\
4	41758.0974441522\\
5	41758.0960693561\\
6	41758.0959514327\\
};


\addlegendentry{TDB-CUR (r=20)}

\end{axis}
\end{tikzpicture}%
    \caption{}
  \end{subfigure}
  \hfill
  \begin{subfigure}[t]{0.23\textwidth}
    \centering
%
%
\definecolor{mycolor1}{rgb}{0.00000,0.44700,0.74100}%
\definecolor{mycolor2}{rgb}{0.85000,0.32500,0.09800}%
\begin{tikzpicture}[scale=0.7]

\begin{axis}[%
width=1.5in,
height=1in,
at={(0.758in,0.481in)},
scale only axis,
unbounded coords=jump,
xmin=1,
xmax=6,
xlabel style={font=\color{white!15!black}},
xlabel={Newton Iteration},
ymin=200,
ymax=300,
yminorticks=true,
ylabel style={font=\color{white!15!black}, align=center},
ylabel={$\sigma_2$},
axis background/.style={fill=white},
xmajorgrids,
ymajorgrids,
yminorgrids,
legend style={
  at={(0,0.355)},
  anchor=north west,
  font=\scriptsize,           
  cells={align=left},           
  align=left,                   
  draw=white!15!black,          
  inner sep=1pt,                
  row sep=1pt,                  
  column sep=5pt,               
}
]
\addplot [color=mycolor1, line width=1.5pt, mark size=1.5pt, mark=x, mark options={solid, mycolor1}, forget plot]
  table[row sep=crcr]{%
1	224.134069249363\\
2	270.499170346304\\
3	272.259418354086\\
4	272.261870881654\\
5	272.261870886435\\
6	272.261870886434\\
};
\addplot [color=mycolor2, dashed, line width=1.5pt, mark size=2.0pt, mark=o, mark options={solid, mycolor2}, forget plot]
  table[row sep=crcr]{%
1	224.134068251013\\
2	270.499166560046\\
3	272.259425692779\\
4	272.261913805441\\
5	272.261874170446\\
6	272.26187105813\\
};



\end{axis}
\end{tikzpicture}%
    \caption{}
  \end{subfigure}
  \begin{subfigure}[t]{0.23\textwidth}
    \centering
%
%
\definecolor{mycolor1}{rgb}{0.00000,0.44700,0.74100}%
\definecolor{mycolor2}{rgb}{0.85000,0.32500,0.09800}%
\begin{tikzpicture}[scale=0.7]

\begin{axis}[%
width=1.5in,
height=1in,
at={(0.758in,0.481in)},
scale only axis,
unbounded coords=jump,
xmin=1,
xmax=6,
xlabel style={font=\color{white!15!black}},
xlabel={Newton Iteration},
ymin=4e-05,
ymax=0.0006,
yminorticks=true,
ylabel style={font=\color{white!15!black}, align=center},
ylabel={$\sigma_{11}$},
axis background/.style={fill=white},
xmajorgrids,
ymajorgrids,
yminorgrids,
legend style={
  at={(0,0.355)},
  anchor=north west,
  font=\scriptsize,           
  cells={align=left},           
  align=left,                   
  draw=white!15!black,          
  inner sep=1pt,                
  row sep=1pt,                  
  column sep=5pt,               
}
]
\addplot [color=mycolor1, line width=1.5pt, mark size=1.5pt, mark=x, mark options={solid, mycolor1}]
  table[row sep=crcr]{%
1	4.15319586957288e-05\\
2	5.2372232539812e-05\\
3	5.31238440634859e-05\\
4	5.31256141161903e-05\\
5	5.31256140913664e-05\\
6	5.31256140604581e-05\\
};
\addplot [color=mycolor2, dashed, line width=1.5pt, mark size=2.0pt, mark=o, mark options={solid, mycolor2}]
  table[row sep=crcr]{%
1	4.15279858254289e-05\\
2	5.23791830015053e-05\\
3	0.000511679834815797\\
4	5.06087938535401e-05\\
5	5.3101926524591e-05\\
6	5.31086109374186e-05\\
};



\end{axis}
\end{tikzpicture}%
    \caption{}
  \end{subfigure}
  \hfill
  \begin{subfigure}[t]{0.23\textwidth}
    \centering
%
%
\definecolor{mycolor1}{rgb}{0.00000,0.44700,0.74100}%
\definecolor{mycolor2}{rgb}{0.85000,0.32500,0.09800}%
\begin{tikzpicture}[scale=0.7]

\begin{axis}[%
width=1.5in,
height=1in,
at={(0.758in,0.481in)},
scale only axis,
unbounded coords=jump,
xmin=1,
xmax=6,
xlabel style={font=\color{white!15!black}},
xlabel={Newton Iteration},
ymin=4.e-06,
ymax=6.e-05,
yminorticks=true,
ylabel style={font=\color{white!15!black}, align=center},
ylabel={$\sigma_{12}$},
axis background/.style={fill=white},
xmajorgrids,
ymajorgrids,
yminorgrids,
legend style={
  at={(0,0.355)},
  anchor=north west,
  font=\scriptsize,           
  cells={align=left},           
  align=left,                   
  draw=white!15!black,          
  inner sep=1pt,                
  row sep=1pt,                  
  column sep=5pt,               
}
]
\addplot [color=mycolor1, line width=1.5pt, mark size=1.5pt, mark=x, mark options={solid, mycolor1}, forget plot]
  table[row sep=crcr]{%
1	4.41084516721044e-06\\
2	7.33483785990528e-06\\
3	7.63249642491242e-06\\
4	7.63330483523477e-06\\
5	7.6333047594194e-06\\
6	7.6333049707319e-06\\
};
\addplot [color=mycolor2, dashed, line width=1.5pt, mark size=2.0pt, mark=o, mark options={solid, mycolor2}, forget plot]
  table[row sep=crcr]{%
1	4.4059021262341e-06\\
2	7.33589617059673e-06\\
3	5.31539336134099e-05\\
4	3.19990168965863e-05\\
5	1.00535758956201e-05\\
6	7.62601691518593e-06\\
};



\end{axis}
\end{tikzpicture}%
    \caption{}
  \end{subfigure}

  \caption{Convergence of singular values for different ranks.}
  \label{fig:pnonlinear2}
\end{figure}
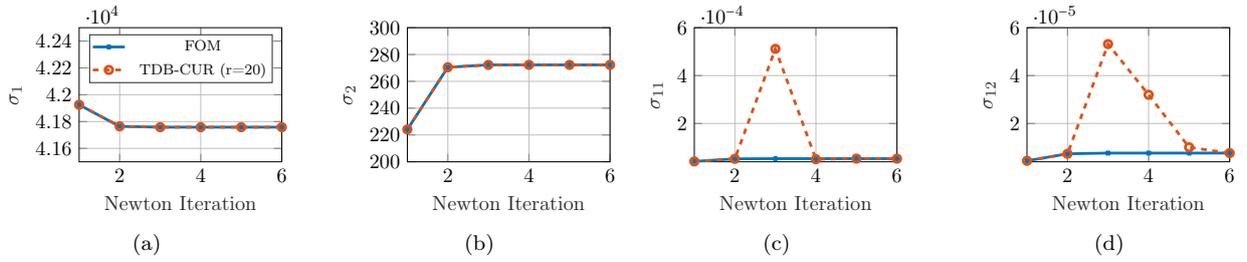

\Cref{fig:pnonlinear2} illustrates the evolution of the first, second, eleventh, and twelfth singular values over the Newton iterations, focused on selected singular values to reveal their convergence behavior in greater detail. The results demonstrate that the proposed TDB–CUR closely follows the FOM for the dominant singular values, particularly \(\sigma_1\) and \(\sigma_2\). These leading singular values have the greatest influence on the overall solution accuracy, whereas higher-order components such as \(\sigma_{11}\) and \(\sigma_{12}\) exhibit minor discrepancies that have negligible effect on the solution accuracy due to their $8$ digit difference from the leading ones.

\subsection{Generalized Lyapunov equation}

We derive the Lyapunov equation for a transient heat conduction problem in two spatial dimensions. The governing equation is
\[
\rho c_p \frac{\partial T}{\partial t} = k \Delta T,
\]
where \( T(x, y, t) \) denotes the temperature, and \( \rho \), \( c_p \), and \( k \) are the density, specific heat capacity, and thermal conductivity, respectively. The equation is discretized in space using the FEM, resulting in a system of linear ordinary differential equations:
\begin{equation}
\bm{M} \bm{\dot{T}} = \bm{K} \bm{T}
\label{eq:heat2Dfem}
\end{equation}
We explain the derivation of the generalized Lyapunov equation for Eq.~\eqref{eq:heat2Dfem} in  \ref{app:lyapunov}.

\subsubsection{Small to moderate size LMEs} \label{subsubsec:case1}
The first case is solving the generalized Lyapunov equation on the domain $\Omega^{(1)}_{(2D)}$ specified by Eq.(\ref{eq:omega}) and shown in Figure~\ref{subig:p2schem2}. 
The material properties are \( k = 50 \,\mathrm{W/(m \cdot K)} \), \( \rho = 780\, \mathrm{kg/m^3} \), and \( c_p = 500\, \mathrm{J/(kg \cdot K)} \).
The details of the generalized Lyapunov equation derivation are provided in \ref{app:lyapunov} and the Eq.~\eqref{lyaFEM1}. The resulting equation is as follows:
\begin{equation} \label{eq:lyaMat}
\bm{A} \bm{X} \bm{B} + \bm{B} \bm{X} \bm{A} + \bm{G} \bm{G}^{\top} = 0,
\end{equation}
where the coefficients $\bm A$ and $\bm B$ are the FEM stiffness and mass matrices, respectively. The matrix $\bm{G}$ is constructed as a rank-1 matrix whose entries correspond to the columnwise sum of the FEM stiffness matrix at all the boundaries. Instead of using the stiffness matrix directly, its columns are summed to reduce the rank, resulting in a rank-1 matrix $\bm{G}\bm{G}^{\top}$. 

We first consider a small-sized problem of $\mathbf X \in \mathbb R^{n \times n}$ with $n=1506$ for which we can compute the FOM solution and perform a detailed comparison of singular value convergence. 
The FOM is solved using GMRES.

To this end, we first solved the full-order generalized Lyapunov equation given above. Since the above equation is not time-dependent, we added a pseudo-time derivative to it as explained in \Cref{sec:ss} and solved the resulting equation using TDB-CUR with the constant $\Delta \tau=10000$  using Eq.~\eqref{eq:tau}.  The TDB-CUR simulation is performed for a fixed-rank of $r=26$. 
\begin{figure}[t]
  \centering
  \begin{subfigure}[t]{0.48\textwidth}
    \centering
    \input{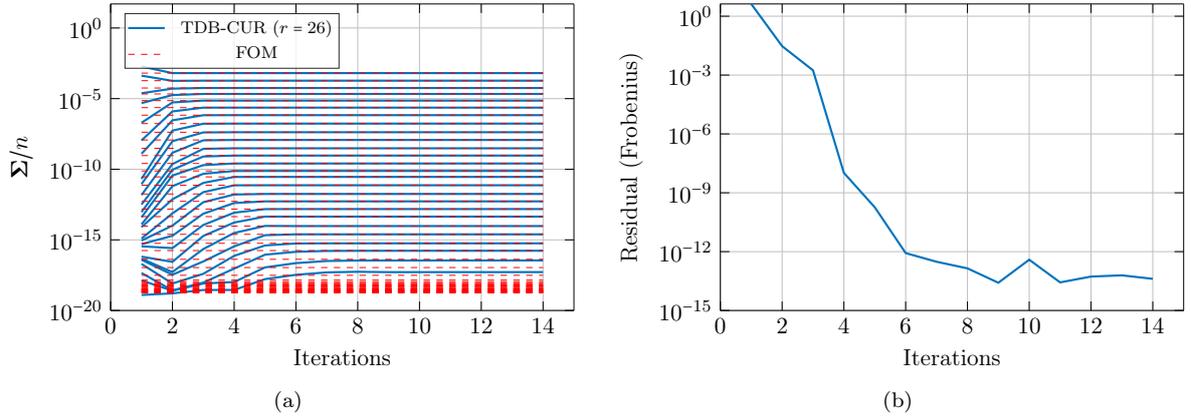}
    \caption{}
  \end{subfigure}
  \begin{subfigure}[t]{0.48\textwidth}
    \centering
%
%
\definecolor{mycolor1}{rgb}{0.00000,0.45000,0.74000}%
\begin{tikzpicture}[scale=0.8]

\begin{axis}[%
width=3in,
height=2.in,
at={(0.874in,0.652in)},
scale only axis,
separate axis lines,
every outer x axis line/.append style={black},
every x tick label/.append style={font=\color{black}},
every x tick/.append style={black},
xmin=0,
xmax=15,
xlabel={Iterations},
every outer y axis line/.append style={black},
every y tick label/.append style={font=\color{black}},
every y tick/.append style={black},
ymode=log,
ymin=1e-15,
ymax=4.20790179549818,
yminorticks=true,
ylabel={Residual (Frobenius)},
axis background/.style={fill=white},
xmajorgrids,
ymajorgrids,
yminorgrids
]
\addplot [color=mycolor1, line width=1.0pt, forget plot]
  table[row sep=crcr]{%
1	4.20790179549818\\
2	0.0297862215355039\\
3	0.00174578483346968\\
4	1.03223751377456e-08\\
5	1.81805268157919e-10\\
6	8.43726661772864e-13\\
7	3.04839872464414e-13\\
8	1.40451806647278e-13\\
9	2.60038088992996e-14\\
10	3.88531999900814e-13\\
11	2.71030865615026e-14\\
12	5.35915575698786e-14\\
13	6.19265521580406e-14\\
14	4.11274051110409e-14\\
};
\end{axis}
\end{tikzpicture}%
    \caption{}
  \end{subfigure}
  \caption{Plate with hole for $n = 1506$. Panel (a) shows the singular values over the number of iterations for TDB-CUR with  ($r=26$) and the first 50 singular values of FOM as horizontal dashed lines. The right panel shows the residual of the TDB-CUR solver.}
  \label{fig:p3sin3}
\end{figure}
\Cref{fig:p3sin3} shows the singular values TDB-CUR method with those of the FOM. The FOM singular values are computed by performing SVD on the solution matrix, while the TDB-CUR singular values are obtained iteratively. It can be observed from Panel~(a) that all singular values after 8 iterations have converged. Panel~(b) shows the decay of the residual in the Frobenius norm over the number of iterations. The residual exhibits rapid convergence, decreasing by several orders of magnitude within the first few iterations and reaching to machine precision. \Cref{fig:p3sin3} confirms that the proposed TDB-CUR method is able to solve the Lyapunov equation of heat conduction effectively.

\begin{figure}[t]
  \centering
  \begin{subfigure}[t]{0.48\textwidth}
    \centering
    \input{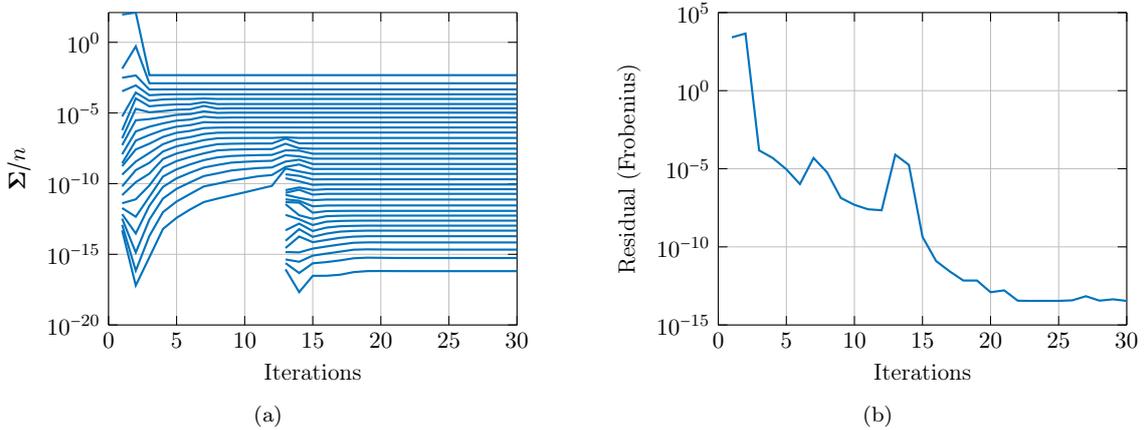}
    \caption{}
  \end{subfigure}
  \begin{subfigure}[t]{0.48\textwidth}
    \centering
%
%
\definecolor{mycolor1}{rgb}{0.00000,0.45000,0.74000}%
\begin{tikzpicture}[scale=0.8]

\begin{axis}[%
width=2.6451in,
height=2.0385in,
at={(0.874in,0.652in)},
scale only axis,
separate axis lines,
every outer x axis line/.append style={black},
every x tick label/.append style={font=\color{black}},
every x tick/.append style={black},
xmin=0,
xmax=30,
xlabel={Iterations},
every outer y axis line/.append style={black},
every y tick label/.append style={font=\color{black}},
every y tick/.append style={black},
ymode=log,
ymin=1e-15,
ymax=100000,
yminorticks=true,
ylabel={Residual (Frobenius)},
axis background/.style={fill=white},
xmajorgrids,
ymajorgrids,
yminorgrids
]
\addplot [color=mycolor1, line width=1.0pt, forget plot]
  table[row sep=crcr]{%
1	2575.57952539956\\
2	4553.16012640134\\
3	0.000150341720719333\\
4	4.91610966843752e-05\\
5	9.24083080603784e-06\\
6	1.04763337122515e-06\\
7	4.91161788789862e-05\\
8	5.78433504733196e-06\\
9	1.35588469183141e-07\\
10	4.91903813234766e-08\\
11	2.54856213600313e-08\\
12	2.20266748405612e-08\\
13	7.96995855918171e-05\\
14	1.76537033077384e-05\\
15	4.38511429024978e-10\\
16	1.23634679679977e-11\\
17	2.67768156441729e-12\\
18	6.98185628815259e-13\\
19	7.02139983450401e-13\\
20	1.24898975044423e-13\\
21	1.64468868219873e-13\\
22	3.56763433429539e-14\\
23	3.46711800574895e-14\\
24	3.50254751346246e-14\\
25	3.48160972307042e-14\\
26	3.77357559246675e-14\\
27	6.93974931582025e-14\\
28	3.60508214145658e-14\\
29	4.39158299763194e-14\\
30	3.45796180260472e-14\\
};
\end{axis}
\end{tikzpicture}%
    \caption{}
  \end{subfigure}
  \caption{Generalized Lyapunov equation solved on the two-dimensional domain $\Omega$ (shown in Figure \ref{subig:p2schem2})  for $n = 102366$ and $r=36$. Panel~(a) shows the singular values over the number of iterations of TDB-CUR. Panel~(b) shows the residual of the TDB-CUR solver over iterations.}
  \label{fig:p3sin4}
\end{figure}

\Cref{fig:p3sin4} presents results for the same problem as the \Cref{fig:p3sin3}, but with a significantly larger problem size. Here, the number of degrees of freedom is $n=102366$. For this problem, computing the FOM solution requires storing $n^2 \approx 10^{10}$ entries, which is not feasible. We employ a rank-adaptive TDB-CUR strategy by initially selecting a rank of  $r=20$ and subsequently increasing it in increments of $\Delta r = 16$. Choosing a very large rank from the outset would impose significant memory requirements due to the size of the Krylov subspace. By contrast, starting with a modest rank and gradually increasing it reduces memory usage, since the converged solution obtained with the initial rank $r=20$ provides an excellent initial guess for the higher-rank GMRES solver. Consequently, the GMRES solver is able to reduce the residual using a much smaller Krylov subspace.  

Table \ref{tab:gmres_results} presents convergence results for different configurations of rank \(r\), the maximum number of bases for the Krylov subspace dimension \(m\), and the maximum memory usage \(n \times r \times m\), showing their impact on residual, iteration count, and wall time under fixed and variable choice of pseudo-time values. Note that each Krylov vector for the rank of $r$ requires storing $n\times r$ entries.  To generate table \ref{tab:gmres_results}, we use a constant pseudo-time value of \( \Delta \tau = 10000 \) and \(n=5000\) is used.

The parameter \( r \times m \) reflects the memory footprint of the proposed reduced-order model. Higher values of \( r \) and \( m \) imply greater memory usage. As shown in the table, increasing \( m \) leads not only to higher memory demand but also to significantly longer wall times. On the other hand, very small values of \( m \) can also increase wall time, resulting in a large number of GMRES restarts. 
The  GMRES iteration reported in the table represents the total number of GMRES iterations and can include multiple restarts.

The solutions with higher ranks are more accurate and generally more costly to solve. In all cases considered, the number of CUR iterations remains relatively low, showing the fast convergence of CUR iterations.

\begin{table}[h]
\centering
\small
\begin{tabular}{|p{1.5cm}|p{1.cm}|p{2 cm}|p{2.0cm}|p{2.0cm}|p{2.0cm}|p{2.0cm}|}
\hline
\textbf{Memory} $(m \times r)$ & \textbf{Rank} ($r$) & \textbf{Krylov Dimension} ($m$) & \textbf{CUR \ \ \ \ \ \ \ Iterations} & \textbf{GMRES Restarts} & \textbf{GMRES Iterations} & \textbf{Wall-time} (seconds)  \\
\hline
60    & 10 & 6   & 6  & 140  & 855   & 2.9   \\
180   & 10 & 18  & 6  & 36   & 670   & 3.0   \\
600   & 10 & 60  & 6  & 6    & 491   & 4.3   \\
6000  & 10 & 600 & 6  & 0    & 460   & 7.4   \\
\hline
60    & 20 & 3   & 10 & 1361 & 4094  & 10.9  \\
180   & 20 & 9   & 10 & 225  & 2060  & 8.1   \\
600   & 20 & 30  & 11 & 46   & 1518  & 15.2  \\
6000  & 20 & 300 & 10 & 0    & 1034  & 38.7  \\
\hline
60    & 30 & 2   & 10 & 5535 & 11079 & 25.7  \\
180   & 30 & 6   & 10 & 674  & 4073  & 17.3  \\
600   & 30 & 20  & 11 & 117  & 2421  & 27.8  \\
6000  & 30 & 200 & 11 & 4    & 1861  & 147.8 \\
\hline
\end{tabular}
\caption{Convergence results for different values of rank ($r$) and the Krylov subspace dimension ($m$). For this study, $\bm X \in \mathbb R^{n \times n}$, where  $n=5000$ and $\Delta \tau=10000$. The first column reports the number of vectors of length $n$ stored in memory. }
\label{tab:gmres_results}
\end{table}

\subsubsection{Large-scale LME} \label{subsubsec:case2}
\begin{figure}
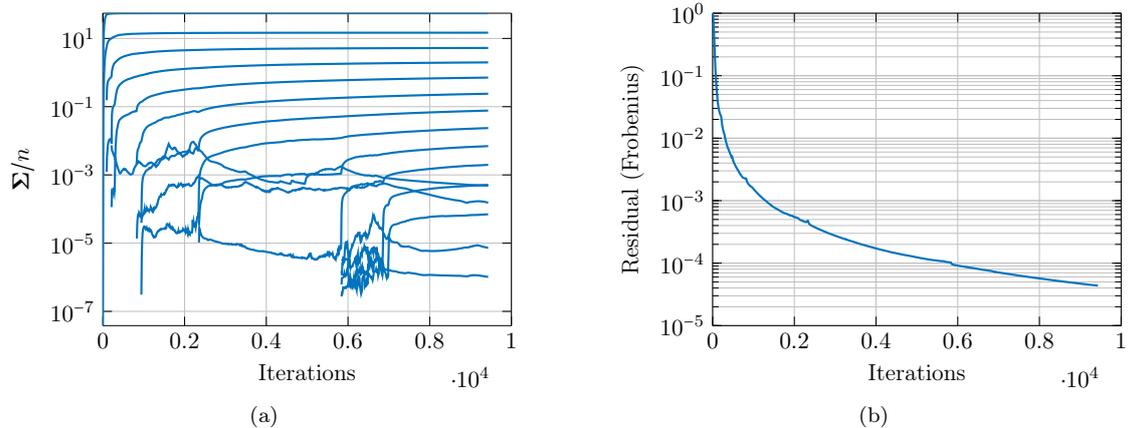

  \centering
  \begin{subfigure}[t]{0.48\textwidth}
    \centering
    \input{FigP3Sin1}
    \caption{ }
  \end{subfigure}
  \begin{subfigure}[t]{0.48\textwidth}
    \centering
    \input{FigP3Res1}
    \caption{ }
  \end{subfigure}

  \caption{Square for $n = 4,242684$. Panel~(a) shows the singular values over the number of iterations for TDB-CUR. Panel~(b) shows the residual of the TDB-CUR solver over iterations.}
  \label{fig:p3sin1_v2}
\end{figure}

The computational domain for the second case is a square.
\begin{equation}
\Omega^{(2)}_{2D} = \left\{ (x,y) \in [-1,1] \times [-1,1] \right\}.
\end{equation}
The material properties are \( k = 1 \,\mathrm{W/(m \cdot K)} \), \( \rho = 1\, \mathrm{kg/m^3} \), and \( c_p = 1\, \mathrm{J/(kg \cdot K)} \). Similar to the Subsection~\ref{subsubsec:case1} the governing matrix equation is
\begin{equation} \label{eq:lyaMat2}
\bm{A} \bm{X} \bm{B} + \bm{B} \bm{X} \bm{A} + \bm{G} \bm{G}^{\top} = 0.
\end{equation}
The known matrices in Eq.~\eqref{eq:lyaMat2} are the same as those in Eq.~\eqref{eq:lyaMat} and are described therein.
\Cref{fig:p3sin1_v2} presents results for the square geometry case, where the spatial resolution is significantly higher with \( n = 4{,}242{,}684 \).
To solve this problem, the solver starts from an initial rank of $r=2$ and increases it in increments of $\Delta r=2$. The solver achieves the target residual tolerance at $r=16$. A variable pseudo-time step is used according to Eq.~\eqref{eq:tau}, with a value of $\Delta \tau_{\text{max}} = 10000$. A preconditioner is not applied; rather, the introduction of a pseudo-time step facilitates smaller GMRES updates at each $\Delta \tau$ increment, though this comes at the expense of increased iteration counts.

This test highlights the scalability of the proposed TDB-CUR method in solving large-scale problems. In the context of solving the Lyapunov equation, this corresponds to a system of approximately \( n^2 \approx 1.8 \times 10^{13} \) algebraic equations.

\section{Conclusion} \label{sec:conc}
In this paper, we present an iterative CUR-based methodology for low-rank approximation and solution of linear matrix equations (LMEs). The CUR decomposition constructs a rank-$r$ approximation of a matrix by sampling $r$ of its columns and rows. When applied to an LME, this decomposition reformulates the large-scale problem into two smaller subproblems: one involving a subset of columns and the other a subset of rows. We propose a fixed-point iterative scheme that removes the dependencies of these subcolumns and subrows on the remaining columns and rows.

At each CUR iteration, two sets of quantities are updated:
(i) the column and row indices, determined adaptively using the discrete empirical interpolation method (DEIM), and
(ii) the coefficient matrices that relate the unsampled columns and rows to the sampled ones.

Numerical experiments demonstrate that the algorithm converges rapidly, typically within 5–10 CUR iterations. The resulting subproblems are solved efficiently using GMRES, allowing the approach to scale to large LMEs.

The proposed framework is agnostic to the number of terms in the LME and can accommodate non-Sylvester equations, such as LMEs involving the Hadamard product between the coefficient and target matrices. We further apply the methodology to implicit time integration of matrix differential equations (MDEs) on low-rank manifolds and to the solution of nonlinear matrix equations via Newton’s method as well as generalized Lyapunov equations. 

An important direction for future research is the design of preconditioners aimed at reducing the number of GMRES iterations and enhancing overall solver efficiency. 

\section*{Acknowledgments}
This work is supported by the Air Force Office of Scientific Research, award no. FA9550-25-1-0039 and award no. FA9550-25-1-0307. 
\vskip2pc

\appendix

\section{DEIM algorithm}
The DEIM pseudocode is presented via Algorithm \ref{alg:DEIM}. This algorithm is adopted from \cite{CS10}. 

\begin{algorithm}
\begingroup
\fontsize{9pt}{9pt}\selectfont
\SetAlgoLined
\KwIn{$\mathbf{U}_{p}=\left[\begin{array}{llll}\mathbf{u}_{1} & \mathbf{u}_{2} & \cdots & \mathbf{u}_{p}\end{array}\right]$}
\KwOut{$\bm I_{p}$}
$\left[\rho, \bm I_{1}\right]=\max \left|\mathbf{u}_{1}\right|$ \hspace{15mm} $\rhd$ choose the first index\;
$\mathbf{P}_{1}=\left[\mathbf{e}_{\bm I_{1}}\right]$ \hspace{24mm} $\rhd$ construct first measurement matrix\;
\For{$i=2$ \KwTo $p$}{
$\mathbf{P}_{i}^{T} \mathbf{U}_{i} \mathbf{c}_{i}=\mathbf{P}_{i}^{T} \mathbf{u}_{i+1}$ \hspace{7.5mm}  $\rhd$ calculate $c_{i}$\;
$\mathbf{R}_{i+1}=\mathbf{u}_{i+1}-\mathbf{U}_{i} \mathbf{c}_{i}$ \hspace{7mm}  $\rhd$ compute residual\;
$\left[\rho, \bm I_{i}\right]=\max \left|\mathbf{R}_{i+1}\right|$ \hspace{7mm} $\rhd$ find index of maximum residual\;
$\mathbf{P}_{i+1}=\left[\begin{array}{ll}\mathbf{P}_{i} & \mathbf{e}_{\bm I_{i}}\end{array}\right]$ \hspace{6mm} $\rhd$ add new column to measurement matrix\;}
\caption{\texttt{DEIM} Algorithm \cite{CS10}}
\label{alg:DEIM}

\endgroup
\end{algorithm}

\section{\label{Stable_CUR}Stable CUR Algorithm for general and symmetric matrices}

The stable CUR pseudocode is presented via
\cref{alg:SCUR}
and we refer to
\cite{DPNFB23}
for more details.
\begin{algorithm}
\begingroup
\fontsize{9pt}{9pt}\selectfont
\SetAlgoLined
\KwIn{$\mathbf{X}(:, \mathbf{q}) \in \mathbb{R}^{n_1 \times r}$, $\mathbf{X}(\mathbf{p}, :) \in \mathbb{R}^{r \times n_2}$, $\mathbf{q}$, $\mathbf{p}$}
\KwOut{$\mathbf{U}$, $\boldsymbol \Sigma$, $\mathbf{Y}$}
$\mathbf{U}_q = \texttt{SVD}(\mathbf{X}(:, \mathbf{q}))$ \hspace{8mm} $\rhd$ Compute the economy SVD of $\mathbf{X}(:, \mathbf{q})$\;
$\mathbf{Y}_p = \texttt{SVD}(\mathbf{X}(\mathbf{p}, :))$ \hspace{6mm} $\rhd$ Compute the economy SVD of $\mathbf{X}(\mathbf{p}, :)$\;
$\mathbf{C} = \mathbf{U}_s(\mathbf{p}, :)^\dagger \, \mathbf{X}(\mathbf{p}, \mathbf{q}) \, (\mathbf{Y}_p(\mathbf{q}, :)^\dagger)^\top$ \hspace{1mm} $\rhd$ Compute $\mathbf{C}$\;
$\mathbf{R}_U, \boldsymbol \Sigma, \mathbf{R}_Y = \texttt{SVD}(\mathbf{C})$ \hspace{11.5mm} $\rhd$ Compute SVD of $\mathbf{C}$\;
$\mathbf{U} = \mathbf{U}_q \mathbf{R}_U$ \hspace{23mm} $\rhd$ In-subspace rotation of $\mathbf{R}_U$\;
$\mathbf{Y} = \mathbf{Y}_p \mathbf{R}_Y$ \hspace{23mm} $\rhd$ In-subspace rotation of $\mathbf{R}_Y$\;
\caption{Stable CUR Algorithm}
\label{alg:SCUR}
\endgroup
\end{algorithm}
\begin{algorithm}
\begingroup
\fontsize{9pt}{9pt}\selectfont
\SetAlgoLined
\KwIn{$\mathbf{U}_R^{k-1} \in \mathbb{R}^{n_1 \times r}$, $\mathbf{Y}_R^{k-1} \in \mathbb{R}^{n_2 \times r}$}
\KwOut{$\mathbf{U}_R^{k}$, $\boldsymbol \Sigma_R^{k}$, $\mathbf{Y}_R^{k}$}
$\bm{q}_R \gets \texttt{DEIM}(\mathbf{Y}_R^{k-1})$ \hspace{7mm} $\rhd$ column indices selection via DEIM\;
$\bm{p}_R \gets \texttt{DEIM}(\mathbf{U}_R^{k-1})$ \hspace{7mm} $\rhd$ row indices selection via DEIM\;
$\mathbf{R}(:, \bm{q}_R) = \mathbf{A}_1 \mathbf{X} \mathbf{B}_1(:, \bm{q}_R) + \mathbf{A}_2 \mathbf{X} \mathbf{B}_2(:, \bm{q}_R) - \mathbf{C}(:, \bm{q}_R)$\;
$\mathbf{R}(\bm{p}_R, :) = \mathbf{A}_1(\bm{p}_R, :) \mathbf{X} \mathbf{B}_1 + \mathbf{A}_2(\bm{p}_R, :) \mathbf{X} \mathbf{B}_2 - \mathbf{C}(\bm{p}_R, :)$\;
$\bm{U}_R^{k}$, $\bm \Sigma_R^{k}$, $\bm{Y}_R^{k}$ $\gets \texttt{Stable\_CUR\_Algorithm}(\bm R(\bm{p}_R, :), \; \bm R(:, \bm{q}_R))$ \hspace{16mm} $\rhd$ Update the residual\;
\caption{Low Rank Approximation of Residual}
\label{alg:res}
\endgroup
\end{algorithm}
\begin{algorithm}
\begingroup
\fontsize{9pt}{9pt}\selectfont
\SetAlgoLined
\KwIn{$\mathbf{U}_\Delta^{k-1} \in \mathbb{R}^{n_1 \times r}$, $\mathbf{Y}_\Delta^{k-1} \in \mathbb{R}^{n_2 \times r}$}
\KwOut{$\mathbf{U}_\Delta^{k}$, $\boldsymbol \Sigma_\Delta^{k}$, $\mathbf{Y}_\Delta^{k}$}
$\bm{q}_\Delta \gets \texttt{DEIM}(\mathbf{Y}_\Delta^{k-1})$ \hspace{7mm} $\rhd$ column indices selection via DEIM\;
$\bm{p}_\Delta \gets \texttt{DEIM}(\mathbf{U}_\Delta^{k-1})$ \hspace{7mm} $\rhd$ row indices selection via DEIM\;
$\Delta\bm X(:, \bm{q}_\Delta) = \bm X^{k}(:, \bm{q}_\Delta) - \bm X^{k-1}(:, \bm{q}_\Delta)$\;
$\Delta\bm X(\bm{p}_\Delta,:) = \bm X^{k}(\bm{p}_\Delta,:) - \bm X^{k-1}(\bm{p}_\Delta,:)$\;
$\bm{U}_\Delta^{k}$, $\bm \Sigma_\Delta^{k}$, $\bm{Y}_\Delta^{k}$ $\gets \texttt{Stable\_CUR\_Algorithm}(\bm \Delta\bm X(\bm{p}_\Delta, :), \; \bm \Delta\bm X(:, \bm{q}_\Delta))$ \hspace{16mm} $\rhd$ Update $\Delta\bm X$\;
\caption{Low Rank Approximation of $\Delta\bm X = \bm X^{k}-\bm X^{k-1}$}
\label{alg:delx}
\endgroup
\end{algorithm}
\section{Finite element method and time integration for the 3D heat equation} \label{app:fem}

In this section, we discretize the three-dimensional heat conduction equation using FEM for the 3D computational domain $\Omega_{3D}$. We assume the solution can be approximated by a separable basis expansion in space and a coefficient function of time:
\begin{equation}
T(x, y, z, t) \approx \sum_{j=1}^{N_z} \sum_{i=1}^{N_{xy}} T_{i,j}(t) \, \phi_i(x, y) \, \psi_j(z),
\end{equation}
where \( \phi_i(x, y) \) are basis functions defined over the two-dimensional domain \( \Omega^{(1)}_{2D} \subset \mathbb{R}^2 \), and \( \psi_j(z) \) are basis functions along the \( z \)-direction. The coefficients \( T_{i,j}(t) \) represent the unknown time-dependent coefficients.

By applying the Galerkin method, we obtain the following semi-discrete system of equations:
\begin{equation} \label{eq:fem2}
\bm M_{xy} \dot{\bm T} \bm M_z = \bm K_{xy} \bm T \bm M_z + \bm M_{xy} \bm T \bm K_z + \bm Q,
\end{equation}
where The coefficients matrices \( \bm{M}_{xy}, \bm{M}_z, \bm{K}_{xy}\), and \(\bm{K}_z \) are defined as:
\begin{equation}
(\bm M_{xy})_{i,k} = \int_{\Omega^{(1)}_{2D}} \phi_i(x,y) \phi_k(x,y) \,dx \, dy,
\end{equation}
\begin{equation}
(\bm M_z)_{j,l} = \int_{0}^{1} \psi_j(z) \psi_l(z) \,dz,
\end{equation}
\begin{equation}
(\bm K_{xy})_{i,k} = \int_{\Omega^{(1)}_{2D}} \nabla \phi_i(x,y) \cdot \nabla \phi_k(x,y) \,dx \, dy,
\end{equation}
\begin{equation}
(\bm K_z)_{j,l} = \int_{0}^{1} \frac{d\psi_j}{dz} \frac{d\psi_l}{dz} \,dz.
\end{equation}
The matrix \( \bm{Q} \) accounts for external forcing terms, given by:
\begin{equation}
(\bm Q)_{i,j} = \int_{\Omega^{(1)}_{2D}} \int_{0}^{1} q(x,y,z) \phi_i(x,y) \psi_j(z) \,dx \,dy \,dz.
\end{equation}
We rewrite Eq.~\eqref{eq:fem2} using the matrices $\bm {A}_i$ and $\bm {B}_i$ in a form similar to the generalized Sylvester equation as follows:
\begin{equation}
\label{syl_3dheat}
\bm A_0  \frac{d \bm X}{d t} \bm B_0 = \bm A_1 \bm X \bm B_1 + \bm A_2 \bm X \bm B_2 + \bm Q,
\end{equation}
where $\bm{A}_0$ and $\bm{A}_2$ correspond to $\bm{M}_{xy}$, $\bm{B}_0$ and $\bm{B}_1$ correspond to $\bm{M}_{z}$, and $\bm{A}_1$ and $\bm{B}_2$ correspond to $\bm{K}_{xy}$ and $\bm{K}_{z}$, respectively.

In the following, we describe the incorporation of Dirichlet boundary conditions into our numerical framework while preserving the symmetry of the coefficient matrices arising from the FEM discretization. The objective is to solve the governing equations only within the interior domain, eliminating equations associated with boundary points, since values of $\bm{X}$ are already known at these locations.
To enforce this for the Dirichlet boundary conditions, we partition the coefficients and unknown matrices into four blocks, corner points ($\bm{A}^{B_{xy}B_{xy}}, \bm{B}^{B_{z}B_{z}}, \text{and } \bm{X}^{B_{xy}B_{z}}$), interior points ($\bm{A}^{i_{xy}i_{xy}}, \bm{B}^{i_{z}i_{z}}, \text{and } \bm{X}^{i_{xy}i_{z}}$), and boundary points excluding corners
(\(\bm{A}^{i_{xy}B_{xy}}, \bm{A}^{B_{xy}i_{xy}}, \bm{B}^{i_{z}B_{z}}, \bm{B}^{B_{z}i_{z}},\) \(\bm{X}^{i_{xy}B_{z}}, \text{and } \bm{X}^{B_{xy}i_{z}}\)),
to distinguish between different points in the grid.
The contributions from boundary grid points appear as source terms on the right-hand side of the system of equations. In the following, we show how to incorporate these boundary terms in the matrix–matrix product and then rewrite the full equation accordingly.
\begin{equation}
\label{eq:matbound}
\bm{A} \bm{X} \bm{B} =
\begin{bmatrix}
\bm{A}^{B_{xy}B_{xy}} & \bm{A}^{B_{xy}i_{xy}} \\
\bm{A}^{i_{xy}B_{xy}} & \bm{A}^{i_{xy}i_{xy}}
\end{bmatrix}
\begin{bmatrix}
\bm{X}^{B_{xy}B_{z}} & \bm{X}^{B_{xy} i_z} \\
\bm{X}^{i_{xy}B_z} & \bm{X}^{i_{xy} i_z}
\end{bmatrix}
\begin{bmatrix}
\bm{B}^{B_{z}B_{z}} & \bm{B}^{B_{z}i_{z}} \\
\bm{B}^{i_{z}B_{z}} & \bm{B}^{i_{z}i_{z}}
\end{bmatrix}
\end{equation}
Since our goal is to solve for $\bm{X}^{i_{xy} i_z}$, we retain only the second block row of Eq.~\eqref{eq:matbound} (interior rows) and disregard the first as follows:
\[
\begin{bmatrix}
\bm{A}^{i_{xy}B_{xy}} & \bm{A}^{i_{xy}i_{xy}}
\end{bmatrix}
\begin{bmatrix}
\bm{X}^{B_{xy} B_z} \\
\bm{X}^{i_{xy} B_z}
\end{bmatrix}
\bm{B}^{B_{z}i_{z}}
+
\begin{bmatrix}
\bm{A}^{i_{xy}B_{xy}} & \bm{A}^{i_{xy}i_{xy}}
\end{bmatrix}
\begin{bmatrix}
\bm{X}^{B_{xy} i_z} \\
\bm{X}^{i_{xy} i_z}
\end{bmatrix}
\bm{B}^{i_{z}i_{z}}
\]
Using the second block row, we express the full equation as follows. 
Note that the left-hand side contains a time-derivative term, which vanishes on Dirichlet boundaries.
\begin{equation}
\label{syl_3dheatODE}
\bm{A}_0^{i_{xy}i_{xy}}
\dot{\bm{X}}^{i_{xy} i_z}
\bm{B}_0^{i_{z}i_{z}}
=
\sum_{j=1}^{2} 
    \left\{
\begin{bmatrix}
\bm{A}_j^{i_{xy}B_{xy}} & \bm{A}_j^{i_{xy}i_{xy}}
\end{bmatrix}
\begin{bmatrix}
\bm{X}^{B_{xy} i_z} \\
\bm{X}^{i_{xy} i_z}
\end{bmatrix}
\bm{B}_j^{i_{z}i_{z}}
+
\begin{bmatrix}
\bm{A}_j^{i_{xy}B_{xy}} & \bm{A}_j^{i_{xy}i_{xy}}
\end{bmatrix}
\begin{bmatrix}
\bm{X}^{B_{xy} B_z} \\
\bm{X}^{i_{xy} B_z}
\end{bmatrix}
\bm{B}_j^{B_{z}i_{z}}
\right\}
 + \bm{Q}
\end{equation}
We employ three different time integrations to discretize the Eq.~\eqref{syl_3dheatODE} on time.
The first-order implicit Euler differentiation formula is given by: 
\begin{equation}
\label{euler}
\bm{A}_0 \dot{\bm{X}} \bm{B}_0 \approx \bm{A}_0 \frac{\bm{X}_{n+1} - \bm{X}_n}{\Delta t} \bm{B}_0.
\end{equation}
The second-order BDF2 formula is given by:
\begin{equation}
\label{bdf2}
\bm{A}_0 \dot{\bm{X}} \bm{B}_0 \approx \bm{A}_0 \frac{3\bm{X}_{n+1} - 4\bm{X}_n + \bm{X}_{n-1}}{2\Delta t} \bm{B}_0.
\end{equation}
The third-order BDF3 formula is given by:
\begin{equation}
\label{bdf3}
\bm{A}_0 \dot{\bm{X}} \bm{B}_0 \approx \bm{A}_0 \frac{11\bm{X}_{n+1} - 18\bm{X}_n + 9\bm{X}_{n-1} - 2\bm{X}_{n-2}}{6\Delta t} \bm{B}_0.
\end{equation}
By substituting Eq.~\eqref{euler}, Eq.~\eqref{bdf2}, and Eq.~\eqref{bdf3} into Eq.~\eqref{syl_3dheatODE}, we obtain three systems of equations, which we solve using the proposed method (see Subsection~\ref{subsec:nearOpt}).

\section{Finite element method for the 3D steady state heat equation with radiation}
\label{app:fem_rad}

The discrete matrix form of the steady-state heat equation with radiation is derived similarly to Eq.~\eqref{eq:fem2}: 
\begin{equation}
\label{eq:fem_rad}
\bm{K}_{xy} \bm{T} \bm{M}_z + \bm{M}_{xy} \bm{T} \bm{K}_z 
- \epsilon \sigma \left[ \bm{G} \circ \left( \bm{T}^4 - \bm{T}_\infty^4 \right) \right] = 0,
\end{equation}
where \( \bm{G} \) is a binary vector that takes the value 1 at the degrees of freedom associated with the radiative boundary \( \Gamma_{\text{rad}} \), and 0 elsewhere. In this derivation, $\bm T^4$ is interpreted elementwise. 
We rewrite Eq.~\eqref{eq:fem_rad} in a form consistent with our earlier notation, using matrices \( \bm{A}_i \) and \( \bm{B}_i \), as:
\begin{equation}
\label{eq:fem_rad_matrix}
\bm{A}_{1} \bm{X} \bm{B}_1 + \bm{A}_{2} \bm{X} \bm{B}_2 
- \epsilon \sigma \left[ \bm{G} \circ \left( \bm{X}^4 - \bm{X}_\infty^4 \right) \right] = 0.
\end{equation}

The corresponding residual form is:
\begin{equation}
\mathcal{R}(\bm{X}) = \bm{A}_{1} \bm{X} \bm{B}_1 + \bm{A}_{2} \bm{X} \bm{B}_2 
- \epsilon \sigma \left[ \bm{G} \circ \left( \bm{X}^4 - \bm{X}_\infty^4 \right).\right]
\end{equation}

We apply Newton linearization to the nonlinear term \(f(\bm{X}) = \bm{X}^4\):
\begin{align}
f(\bm{X}^{k+1}) &= f(\bm{X}^k + \delta \bm{X}) \approx f(\bm{X}^k) + f'(\bm{X}^k)\delta \bm{X} \notag \\
&= (\bm{X}^k)^4 + 4(\bm{X}^k)^3 \circ \delta \bm{X}.
\end{align}

We linearize the nonlinear equation:
\begin{equation}
\bm{A}_1 \bm{X}^{k+1} \bm{B}_1 + \bm{A}_2 \bm{X}^{k+1} \bm{B}_2 
- \epsilon \sigma \, \bm{G} \circ \left( (\bm{X}^{k+1})^4 - \bm{X}_\infty^4 \right) = 0.
\end{equation}

Using the Newton update \( \bm{X}^{k+1} = \bm{X}^k + \delta \bm{X} \), we expand:
\begin{equation}
\bm{A}_1 (\bm{X}^k + \delta \bm{X}) \bm{B}_1 + \bm{A}_2 (\bm{X}^k + \delta \bm{X}) \bm{B}_2 
- \epsilon \sigma \, \bm{G} \circ \left( (\bm{X}^k)^4 + 4(\bm{X}^k)^3 \circ \delta \bm{X} - \bm{X}_\infty^4 \right) = 0.
\end{equation}

Grouping unknown and known terms gives:
\begin{equation}
\underbrace{\bm{A}_1 \bm{X}^k \bm{B}_1 + \bm{A}_2 \bm{X}^k \bm{B}_2 
- \epsilon \sigma \, \bm{G} \circ \left( (\bm{X}^k)^4 - \bm{X}_\infty^4 \right)}_{\mathcal{R}(\bm{X}^k)} 
+ \bm{A}_1 \delta \bm{X} \bm{B}_1 + \bm{A}_2 \delta \bm{X} \bm{B}_2 
- 4 \epsilon \sigma \, \bm{G} \circ \left((\bm{X}^k)^3 \circ \delta \bm{X}\right) = 0.
\end{equation}

This leads to the Newton step:
\begin{equation}
\label{eq:newton1}
\bm{A}_1 \delta \bm{X} \bm{B}_1 + \bm{A}_2 \delta \bm{X} \bm{B}_2 
+ 4 \delta \epsilon \, \bm{G} \circ \left((\bm{X}^k)^3 \circ \delta \bm{X}\right)
= -\mathcal{R}(\bm{X}^k).
\end{equation}

Similar to ~\ref{app:fem}, we use Eq.~\eqref{eq:matbound} to separate the Dirichlet boundary conditions from Eq.~\eqref{eq:newton1}, since they are given and do not need to be solved. The resulting system is a linear matrix equation, which we then solve using our proposed method.

\section{Generalized Lyapunov Equation for Finite Element Discretization} \label{app:lyapunov}

The derivation begins with the FEM approximation of the transient heat conduction. In particular, in this work, we used standard P1 Lagrange finite elements. This results in a system of equations describing the time evolution of the temperature field:
\begin{equation}
\bm{M} \bm{\dot{T}} = \bm{K} \bm{T},
\end{equation}
where $\bm{K}$ is the stiffness matrix, $\bm{M}$ is the mass matrix, and $\bm{T}$ is the vector containing the values of the degrees of freedom, which, in the present case, represents the temperature at the mesh points. Since the temperature values at the boundary nodes are known at the boundaries with the Dirichlet condition, we partition the system into boundary, interior, and their interaction blocks as follows:
\begin{equation}
\label{FEMmatrixEq}
\begin{bmatrix}
\bm{M}_{BB} & \bm{M}_{Bi} \\
\bm{M}_{iB} & \bm{M}_{ii}
\end{bmatrix}
\begin{bmatrix}
\bm 0 \\
\bm{\dot{T}}_i
\end{bmatrix}
=
\begin{bmatrix}
\bm{K}_{BB} & \bm{K}_{Bi} \\
\bm{K}_{iB} & \bm{K}_{ii}
\end{bmatrix}
\begin{bmatrix}
\bm{T}_B \\
\bm{T}_i
\end{bmatrix}.
\end{equation}
Utilizing the second row of the matrix Eq.~\eqref{FEMmatrixEq} leads to an explicit equation for the temperature at the interior points
\begin{equation}
\label{interiourEq}
\bm{M}_{ii} \bm{\dot{T}}_i = \bm{K}_{iB} \bm{T}_B + \bm{K}_{ii} \bm{T}_i.
\end{equation}
We multiply the Eq.~\eqref{interiourEq} to $\bm{M}_{ii}^{-1}$ to form an equation similar to state space form
\begin{equation}
\bm{\dot{T}}_i = \bm{M}_{ii}^{-1} \bm{K}_{iB} \bm{T}_B + \bm{M}_{ii}^{-1} \bm{K}_{ii} \bm{T}_i.
\end{equation}
The state space and the related Lyapunov equation are
\begin{equation}
\label{eq:lyaCont}
\bm{\dot{x}} = \bm{A} \bm{x} + \bm{B} \bm{u} \quad \longrightarrow \quad \bm{A} \bm{X} + \bm{X} \bm{A}^{\top} + \bm{B} \bm{B}^{\top} = 0,
\end{equation}
where we define $\bm{A}$ and $\bm{B}$, the coefficient of the state space form from our FEM equation
\begin{equation}
\bm{A} = \bm{M}_{ii}^{-1} \bm{K}_{ii}, \quad \bm{B} = \bm{M}_{ii}^{-1} \bm{K}_{iB}.
\end{equation}
By substituting these expressions, we obtain the Lyapunov equation:
\begin{equation}
\label{lyaFEM1}
\bm{M}_{ii}^{-1} \bm{K}_{ii} \bm{X} + \bm{X} (\bm{M}_{ii}^{-1} \bm{K}_{ii})^{\top} + \bm{M}_{ii}^{-1} \bm{K}_{iB} (\bm{M}_{ii}^{-1} \bm{K}_{iB})^{\top} = 0.
\end{equation}
By multiplying the $\bm{M}_{ii}^{-1}$ from left and right, the Eq.~\eqref{lyaFEM1} can be rewritten as:
\begin{equation}
\bm{K}_{ii} \bm{X} \bm{M}_{ii} + \bm{M}_{ii} \bm{X} \bm{K}_{ii} + \bm{K}_{iB} \bm{K}_{iB}^{\top} = 0.
\end{equation}
Since the Lyapunov equation is steady-state, solving it needs a high-dimensional Krylov space, leading to high storage usage. To reduce the required storage, we introduce a pseudo-time derivative and limit the Krylov space to a maximum dimension and restart it to control the storage usage:
\begin{equation}
\dot{\bm  X} = \bm{K}_{ii} \bm{X} \bm{M}_{ii} + \bm{M}_{ii} \bm{X} \bm{K}_{ii} + \bm{K}_{iB} \bm{K}_{iB}^{\top}.
\end{equation}
In practice, retaining the full matrix \(\bm{K}_{iB}\) leads to a high-rank source term \(\bm{B}\bm{B}^\top\), since \(\bm{K}_{iB}\) has a column for each boundary node. This means the degrees of freedom of controlling the system are equal to the grid points at the boundaries, significantly increasing the effective rank of the Lyapunov equation. To reduce the rank of the system, we replace \(\bm{K}_{iB}\), we sum of its columns to construct a vector \(\bm{B}\). This reduces the source term to rank 1, enabling low-rank solvers to operate more efficiently. Physically, this simplification corresponds to coupling the system to a single scalar control input that modulates the aggregate influence of the boundary nodes.

After Euler time integration of the Lyapunov equation with the pseudo-time derivative, we have the following equation:
\begin{equation}
\bm X^{k+1} - \tau \big [\bm{K}_{ii} \bm X^{k+1} \bm{M}_{ii} +  \bm{M}_{ii} \bm X^{k+1} \bm{K}_{ii} \big ] = \bm X^{k} + \tau \bm{B} \bm{B}^{\top}.
\end{equation}
The above equation has the form of \ref{eq:trgse}.
\bibliographystyle{elsarticle-num}
\bibliography{References,Hessam}
\end{document}